\documentclass[11pt, twocolumn]{article}

\usepackage{scicite}
\usepackage{mathptmx}
\usepackage[top=0.69in, bottom=0.56in, left=.56in, right=.56in]{geometry}
\usepackage{setspace}

\usepackage{graphicx}
\usepackage{amsmath}
\usepackage{amssymb}
\usepackage{amsfonts}
\usepackage{float}
\usepackage[title]{appendix}

\usepackage{array}
\usepackage{makecell}
\setcellgapes{10pt}

\usepackage[normalem]{ulem}
\usepackage{xcolor}

\def\bx{\mathbf{x}}
\newcommand{\Gv}{\mathbf{G}}
\newcommand{\Nv}{\mathbf{N}}
\newcommand{\Lv}{\mathbf{L}}
\newcommand{\uv}{\mathbf{u}}
\newcommand{\Fv}{\mathbf{F}}
\newcommand{\vv}{\mathbf{v}}
\newcommand{\fv}{\mathbf{f}}
\newcommand{\xv}{\mathbf{x}}
\newcommand{\Bv}{\mathbf{B}}
\newcommand{\bpsi}{\boldsymbol{\psi}}
\newcommand{\bphi}{\boldsymbol{\phi}}

\graphicspath{ {./figures/} }

\title{\textbf{DeepGreen: Deep Learning of Green's Functions\\ for Nonlinear Boundary Value Problems}}
\author{}

\author
{Craig R. Gin$^{1}$\footnote{Co-first authors.}~, Daniel E. Shea$^{2\ast}$, Steven L. Brunton$^{3}$, and J. Nathan Kutz$^{4}$ 

\\
\normalsize{$^{1}$Department of Population Health and Pathobiology, North Carolina State University,}\\
\normalsize{Raleigh, NC 27695, United States}\\
\normalsize{$^{2}$Department of Materials Science and Engineering, University of Washington,}\\
\normalsize{Seattle, WA 98195, United States}\\
\normalsize{$^{3}$Department of Mechanical Engineering, University of Washington,}\\
\normalsize{Seattle, WA 98195, United States}\\
\normalsize{$^{4}$Department of Applied Mathematics, University of Washington,}\\
\normalsize{Seattle, WA 98195, United States}\\
\\
\normalsize{$^\ast$Co-first authors and corresponding authors. Email: crgin@ncsu.edu and sheadan@uw.edu.}
}

\date{}

\begin{document}

\twocolumn[
\begin{@twocolumnfalse}
	\maketitle
	\vspace{-.3in}
	
	\begin{abstract}
		{\bf \noindent 	
Boundary value problems (BVPs) play a central role in the mathematical analysis of constrained physical systems subjected to external forces. 
Consequently, BVPs frequently emerge in nearly every engineering discipline and span problem domains including fluid mechanics, electromagnetics, quantum mechanics, and elasticity.
The fundamental solution, or Green's function, is a leading method for solving linear BVPs that enables facile computation of new solutions to systems under any external forcing. 
However, fundamental Green's function solutions for nonlinear BVPs are not feasible since linear superposition no longer holds. 
In this work, we propose a flexible deep learning approach to solve nonlinear BVPs using a dual-autoencoder architecture.
The autoencoders discover an invertible coordinate transform that linearizes the nonlinear BVP and identifies both a linear operator $L$ and Green's function $G$ which can be used to solve new  nonlinear BVPs.
We find that the method succeeds on a variety of nonlinear systems including nonlinear Helmholtz and Sturm--Liouville problems, nonlinear elasticity, and a 2D nonlinear Poisson equation.
The method merges the strengths of the universal approximation capabilities of deep learning with the physics knowledge of Green's functions to yield a flexible tool for identifying fundamental solutions to a variety of nonlinear systems.\\

		}
	\end{abstract}
\end{@twocolumnfalse}
]

\section{Introduction}

Boundary value problems (BVPs) are ubiquitous in the sciences~\cite{stakgold2000boundary}.   From elasticity to quantum electronics, BVPs have been fundamental in the development and engineering design of numerous transformative technologies of the 20th century.  
Historically, the formulation of many canonical problems in physics and engineering result in {\em linear} BVPs:  from Fourier formulating the heat equation in 1822~\cite{fourier_theorie_1822} to more modern applications such as designing chip architectures in the semi-conductor industry~\cite{jackson,yariv}.  Much of our theoretical understanding of BVPs comes from the construction of the fundamental solution of the BVP, commonly known as the Green's function~\cite{stakgold2011green}.  The Green's function solution relies on a common property of many BVPs:  {\em linearity}.  Specifically, general solutions rely on linear superposition to hold, thus limiting their usefulness in many modern applications where BVPs are often heterogeneous and nonlinear. 
By leveraging modern deep learning, we are able to learn linearizing transformations of BVPs that render {\em nonlinear BVPs linear} so that we can construct the Green's function solution.  Our deep learning of Green's functions, {\em DeepGreen}, provides a transformative architecture for modern solutions of nonlinear BVPs.    

DeepGreen is inspired by recent works which use deep neural networks (DNNs) to discover advantageous coordinate transformations for dynamical systems~\cite{lusch2018deep,championDatadrivenDiscoveryCoordinates2019,Wehmeyer2017,Mardt2018,Takeishi2017nips,Yeung2019,Otto2019,Li2017chaos,dsilva2018parsimonious, gin2020}.  
The universal approximation properties of DNNs~\cite{Cybenko.1989,Hornik.1990} are ideal for learning  coordinate transformations that linearize nonlinear BVPs, ODEs and PDEs.  Specifically, such linearizing transforms fall broadly under the umbrella of Koopman operator theory~\cite{Koopman1931pnas}, which has a modern interpretation in terms of dynamical systems theory~\cite{Mezic2004,Mezic2005nd,budivsic2012applied,Mezic2013arfm}.  There are only limited cases in which Koopman operators can only be constructed explicitly~\cite{brunton_koopman_2016}. However {\em Dynamic Mode Decomposition} (DMD)~\cite{Schmid2010jfm} provides a numerical algorithm for approximating the Koopman operator~\cite{Rowley2009jfm}, with many recent extensions that improve on the DMD approximation~\cite{Kutz2016book}.  More recently, neural networks have been used to construct Koopman embeddings~\cite{lusch2018deep,Wehmeyer2017,Mardt2018,Takeishi2017nips,Yeung2019,Otto2019,Li2017chaos,gin2020}.  This is an alternative to enriching the observables of DMD~\cite{noe2013variational,nuske2014variational,Williams2015jnls,Williams2015jcd,klus2017data,kutzPDE,page2018koopman}.  Thus, neural networks have emerged as a highly effective mathematical tool for approximating complex data~\cite{ml3,GoodfellowDL} with a {\em linear} model.  DNNs have been used in this context to discover time-stepping algorithms for complex systems~\cite{rico1995nonlinear,gonzalez1998identification,rudy2019deep,lange2020fourier,liu2020hierarchical}.  Moreover, DNNs have been used to approximate constitutive models of  BVPs~\cite{huangLearningConstitutiveRelations2020}.

DeepGreen leverages the success of DNNs for dynamical systems to discover coordinate transformations that linearize nonlinear BVPs so that the Green's function solution can be recovered.  This allows for the discovery of the fundamental solutions for nonlinear BVPs, opening many opportunities for the engineering and physical sciences.  DeepGreen exploits physics-informed learning by using autoenconders (AEs) to take data from the original high-dimensional input space to the new coordinates at the intrinsic rank of the underlying physics~\cite{lusch2018deep, pan2019, championDatadrivenDiscoveryCoordinates2019}. 
The architecture also leverages the success of {\em Deep Residual Networks} (DRN)~\cite{he2016deep} which enables our approach to efficiently handle near-identity coordinate transformations~\cite{gin2020}. 

The Green's function constructs the solution to a BVP for any given forcing by superposition.  Specifically, consider the classical linear BVP~\cite{stakgold2011green} 
\begin{equation}\label{luf}
\begin{array}{ll}
    L[v({\bx})] = f({\bx}) 
\end{array}
\end{equation}
where $L$ is a linear differential operator, $f$ is a forcing, ${\bx}\in {\Omega}$ is the spatial coordinate, and ${\Omega}$ is an open set.  The boundary conditions $Bv({\bx})=0$ are imposed on  $\partial {\Omega}$ with a linear operator $B$.  The fundamental solution is constructed by considering the adjoint equation
\begin{equation}
\label{g}
\begin{array}{ll}
    L^\dag [G({\bx,\boldsymbol{\xi}})] = \delta ({\bx}-{\boldsymbol{\xi}})
\end{array}
\end{equation}
where $L^\dag$ is the adjoint operator (along with its associated boundary conditions) and $ \delta ({\bx}-{\bf \xi})$ is the Dirac delta function.  Taking the inner product of (\ref{luf}) with respect to the Green's function gives the fundamental solution
\begin{equation}
\label{eq:green}
    v({\bf x}) = (f({\boldsymbol{\xi}}), G(\boldsymbol{\xi},{\bf x}) ) =
    \int\limits_{{\Omega}} G(\mathbf{\xi},{\bf x}) f(\boldsymbol{\xi}) d\boldsymbol{\xi},
\end{equation}
which is valid for any forcing $f({\bx})$.  Thus once the Green's function is computed, the solution for arbitrary forcing functions can be easily extracted from integration.  This integration represents a superposition of a continuum of delta function forcings that are used to represent $f({\bx})$.

In many modern applications, nonlinearity plays a fundamental role so that the BVP is of the form
\begin{equation}
\begin{array}{ll}
    N[u({\bx})] = F({\bx}) 
\end{array}
\end{equation}
where $N[\cdot]$ is a nonlinear differential operator.  For this case, the principle of linear superposition no longer holds and the notion of a fundamental solution is lost.  However, modern deep learning algorithms allow us the flexibility of learning a coordinate transformation (and their inverses) of the form
\begin{subequations} \label{eqn:coordinate_transform}
\begin{align}
    v &= \bpsi(u), \\
    f &= \bphi(F),
\end{align}
\end{subequations}
such that $v$ and $f$ satisfy the linear BVP (\ref{luf}) for which we generated the fundamental solution (\ref{eq:green}).  This gives a nonlinear fundamental solution through use of this deep learning transformation.

DeepGreen is a {\em supervised learning} algorithm which is ultimately a high-dimensional interpolation problem~\cite{mallat2016understanding} for learning the coordinate transformations $\bpsi(u)$ and $\bphi(F)$.  DeepGreen is enabled by a physics-informed deep autoencoder coordinate transformation which establishes superposition for nonlinear BVPs, thus enabling a Koopman BVP framework.  The learned Green's function enables accurate construction of solutions with new forcing functions in the same way as a linear BVP.  We demonstrate the DeepGreen method on a variety of nonlinear boundary value problems, including a nonlinear 2D Poisson problem, showing that such an architecture can be used in many modern and diverse applications in aerospace, electromagnetics, elasticity, materials, and chemical reactors.

\begin{figure*}[t]
    \centering
    \includegraphics[width=15cm]{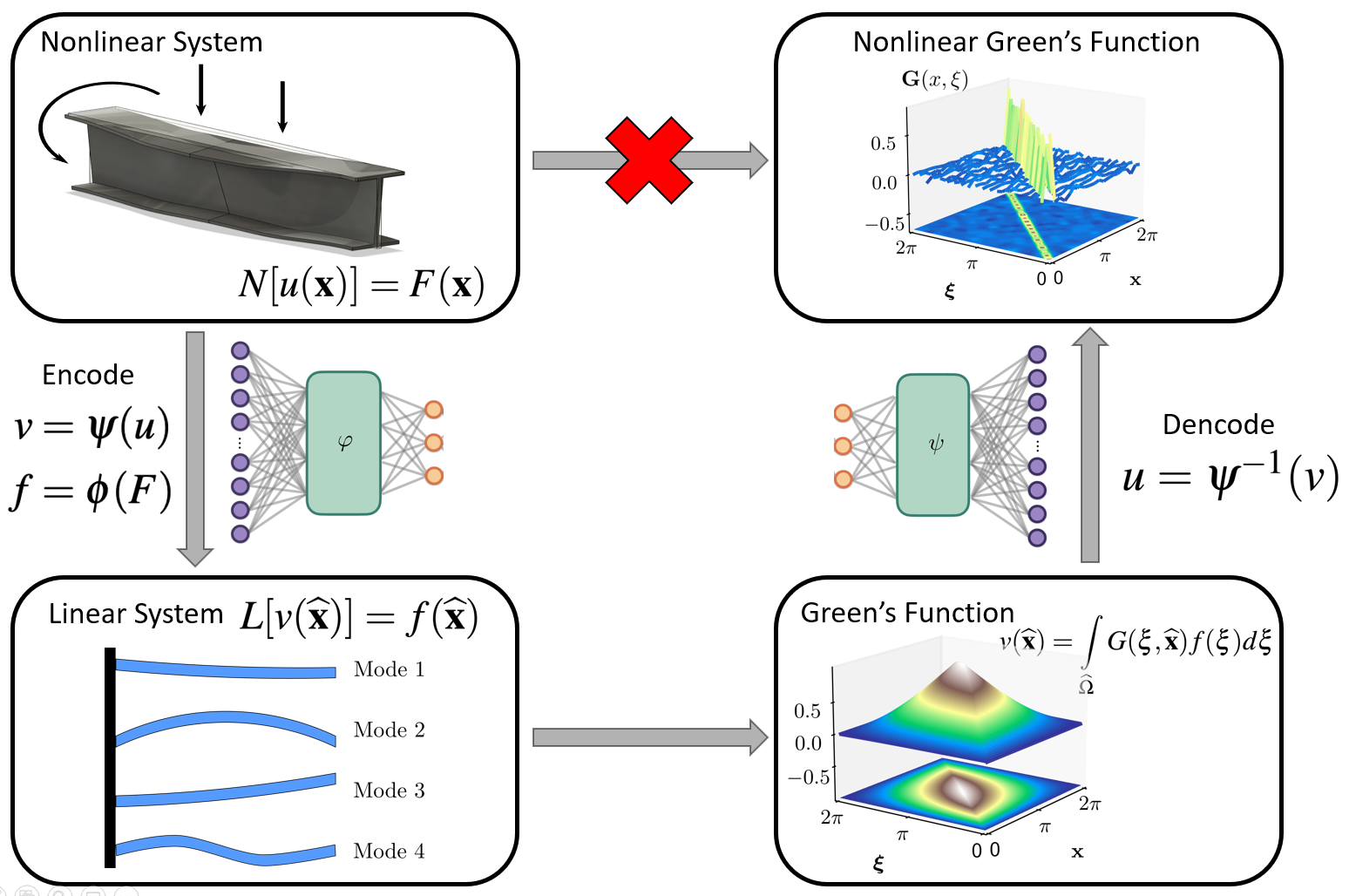}
    \caption{DeepGreen solves nonlinear BVPs by identifying the Green's Function of the nonlinear problem using a deep learning approach with a dual autoencoder architecture. An nonhomogenous linear BVP can be solved using the Green's function approach, but a nonlinear BVP cannot. DeepGreen transforms a nonlinear BVP to a linear BVP, solves the linearized BVP, and then inverse transforms the linear solution to solve the nonlinear BVP.}
    \label{fig:graphical_abstract}
\end{figure*}

\section{Deep Autoencoders for Linearizing BVPs} \label{sec:dae-bvp}

Deep AEs have been used to linearize dynamical systems, which are initial value problems. We extend this idea to BVPs. To be precise, we consider BVPs of the form
\begin{subequations}\label{BVP_continuous}
\begin{align}
    N[u(\xv)] &= F(\xv), & 
    \bx &\in \Omega, \\
    B[u(\xv)] &= 0, & 
    \bx &\in \partial \Omega,
\end{align}
\end{subequations}
where $\Omega$ is a simply connected open set in $\mathbb{R}^n$ with boundary $\partial \Omega$, $N$ is a nonlinear differential operator, $F(\xv)$ is the nonhomogeneous forcing function, $B$ is a boundary condition, and $u(\xv)$ is the solution to the BVP. We wish to find a pair of coordinate transformations of the form \eqref{eqn:coordinate_transform}
%
%
such that $v$ and $f$ satisfy a linear BVP
\begin{subequations}\label{BVP_linear_cont}
\begin{align}
    L[v(\widehat{\bx})] &= f(\widehat{\bx}), & 
    \widehat{\bx} &\in \widehat{\Omega},  \\
    \widehat{B}[v(\widehat{\bx})] &= 0, & 
    \widehat{\bx} &\in \partial \widehat{\Omega},
\end{align}
\end{subequations}
where $L$ is a linear differential operator, $\widehat{\bx}$ is the spatial coordinate in the transformed domain $\widehat{\Omega}$ with boundary $\partial \widehat{\Omega}$.  Because $L$ is linear, there is a Green's function $G(\widehat{\bx},\boldsymbol{\xi})$ such that the solution $v$ to the BVP \eqref{BVP_linear_cont} can be obtained through convolution of the Green's function and transformed forcing function
\begin{equation}\label{Greens}
    v(\widehat{\bx}) = \int\limits_{\widehat{\Omega}} G(\mathbf{\xi},\widehat{\bx}) f(\boldsymbol{\xi}) d\boldsymbol{\xi}.
\end{equation}
The coordinate transformation along with the Green's function of the linearized BVP provide the analog of a Green's function for the nonlinear BVP \eqref{BVP_continuous}. In particular, for a forcing function $F(\xv)$, the transformed forcing function is $f = \bphi(F)$. The solution to the linearized BVP can be obtained using the Green's function $v = \int G(\boldsymbol{\xi},\widehat{\bx}) f(\boldsymbol{\xi}) d\boldsymbol{\xi}$. Then the solution to the nonlinear BVP \eqref{BVP_continuous} is obtained by inverting the coordinate transformation $u=\bpsi^{-1}(v)$  to obtain the solution to the nonlinear BVP, $u(\xv)$.

The question that remains is how to discover the appropriate coordinate transformations $\bpsi$ and $\bphi$. We leverage the universal approximation properties of neural networks in order to learn these transformations. In order to use neural networks, we first need to discretize the BVP. Let $\uv$ be a spatial discretization of $u(\xv)$ and $\Fv$ be a discretization of $F(\xv)$. Then the discretized version of the BVP \eqref{BVP_continuous} is
\begin{subequations}\label{BVP}
\begin{align}
    \Nv[\uv] &= \Fv,  \\
    \Bv[\uv] &= \mathbf{0}.
\end{align}
\end{subequations}
Neural networks $\bpsi_u$ and $\bphi_F$ are used to transform $\uv$ and $\Fv$ to the latent space vectors $\vv$ and $\fv$
\begin{subequations}
\begin{align}
    \vv &= \bpsi_u(\uv), \label{u-transform}\\ 
    \fv &= \bphi_F(\Fv), \label{F-transform}
\end{align}
\end{subequations}
where $\vv$ and $\fv$ satisfy the linear equation
\begin{equation}\label{BVP_linear}
    \Lv \vv = \fv,
\end{equation}
for some matrix $\Lv$, which is also learned. In order to learn invertible transforms $\bpsi_u$ and $\bphi_F$, we construct the problem as a pair of autoencoder networks. 

\begin{figure*}[t]
    \centering
    \includegraphics[width=14cm]{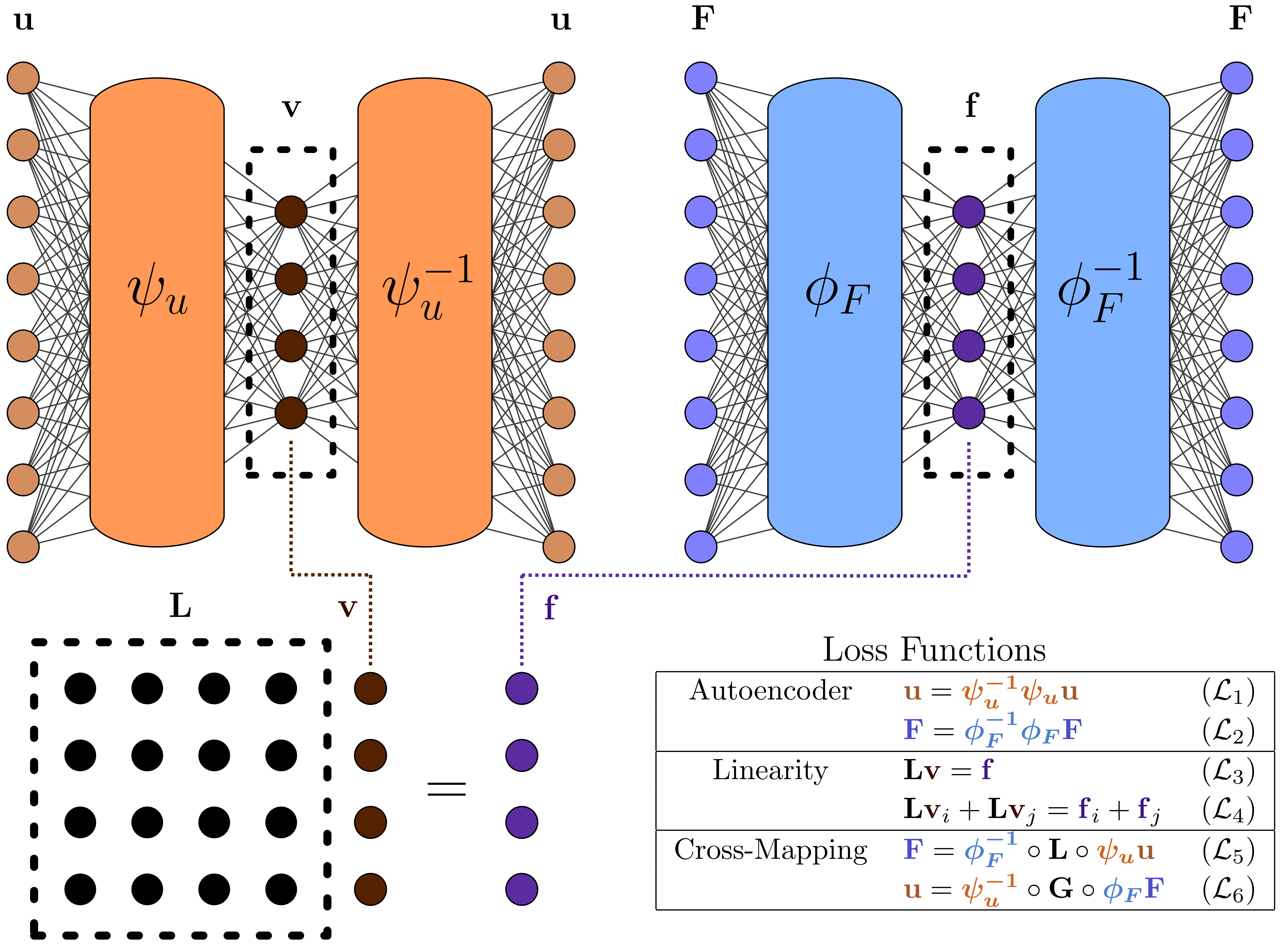}
    \caption{DeepGreen architecture. Two autoencoders learn invertible coordinate transformations that linearize a nonlinear boundary value problem. The latent space is constrained to exhibit properties of a linear system, including linear superposition, which enables discovery of a Green's function for nonlinear boundary value problems.}
    \label{fig:greennet_arch}
\end{figure*}

In this construction, the transforms $\bpsi_u$ and $\bphi_F$ are the encoders and the transform inverses are the decoders. The network architecture and loss functions are shown in Figure \ref{fig:greennet_arch}.
The neural network is trained using numerous and diverse solutions to the nonlinear BVP \eqref{BVP}, which can be obtained with many different forcings $\Fv_k$. Consider a dataset comprised of pairs of discretized solutions and forcing functions $\{\uv_k,\Fv_k\}_{k=1}^N$. The loss function for training the network is the sum of six losses, each of which enforces a desired condition. The loss functions can be split into three categories:
\begin{enumerate}
    \item \textbf{Autoencoder losses:}
    We wish to learn invertible coordinate transformations given by equations \eqref{u-transform} and \eqref{F-transform}. In order to do so, we use two autoencoders. The autoencoder for $\uv$ consists of an encoder $\bpsi_u$ which performs the transformation \eqref{u-transform} and a decoder $\bpsi_u^{-1}$ which inverts the transformation. In order to enforce that the encoder and decoder are inverses, we use the autoencoder loss
    \begin{equation}
         \mathcal{L}_1 = \frac{1}{N} \sum_{k=1}^N \frac{\left\lVert \uv_k - \bpsi_u^{-1} \circ \bpsi_u(\uv_k) \right\rVert_2^2}{\left\lVert \uv_k\right\rVert_2^2}.
    \end{equation}
    Similarly, there is an autoencoder for $\Fv$ where the encoder $\bphi_F$ performs the transformation \eqref{F-transform}. This  transformation also has an inverse enforced by the associated autoencoder loss function
    \begin{equation}
         \mathcal{L}_2 = \frac{1}{N} \sum_{k=1}^N \frac{\left\lVert \Fv_k - \bphi_F^{-1} \circ \bphi_F(\Fv_k) \right\rVert_2^2}{\left\lVert \Fv_k\right\rVert_2^2}.
    \end{equation}
    \item \textbf{Linearity losses:}
    In the transformed coordinate system, we wish for the BVP to be linear so that the operator can be represented by a matrix $\Lv$. The matrix $\Lv$ and the encoded vectors $\vv$ and $\fv$ should satisfy equation \eqref{BVP_linear}. This is enforced with the linear operator loss
    \begin{equation}
         \mathcal{L}_3 = \frac{1}{N} \sum_{k=1}^N \frac{\left\lVert \fv_k - \Lv \vv_k \right\rVert_2^2}{\left\lVert \fv_k\right\rVert_2^2}.
    \end{equation}
    The major advantage of working with a linear operator is that linear superposition holds. We use a linear superposition loss in order to further enforce the linearity of the operator in the latent space
    \begin{equation}
         \mathcal{L}_4 = \frac{1}{N^2} \sum_{j=1}^N \sum_{i=1}^N \frac{\left\lVert (\fv_i+\fv_j) - \Lv (\vv_i+\vv_j) \right\rVert_2^2}{\left\lVert \fv_i + \fv_j \right\rVert_2^2}.
    \end{equation}
    \item \textbf{Cross-mapping losses:}
    The losses described above are theoretically sufficient to find coordinate transformations for $\uv$ and $\Fv$ as well as a linear operator $\Lv$. However, in practice the two autoencoders were not capable of generating the Green's function solution. To rectify this, we add two ``cross-mapping" loss functions that incorporate parts of both autoencoders. The first cross-mapping loss enforces the following mapping from $\uv$ to $\Fv$. First, one of the solutions from the dataset $\uv_k$ is encoded with $\bpsi_u$. This is an approximation for $\vv_k$. This is then multiplied by the matrix $\Lv$, giving an approximation of $\fv_k$. Then the result is decoded with $\bphi_F^{-1}$. This gives an approximation of $\Fv_k$. The $\uv$ to $\Fv$ cross-mapping loss is given by the formula
    \begin{equation}
         \mathcal{L}_5 = \frac{1}{N} \sum_{k=1}^N \frac{\left\lVert \Fv_k - \bphi_F^{-1} \circ \Lv \circ \bpsi_u(\uv_k) \right\rVert_2^2}{\left\lVert \Fv_k\right\rVert_2^2}.
    \end{equation}
    We can similarly define a cross-mapping from $\Fv$ to $\uv$. For a forcing function $\Fv_k$ from the dataset, it is encoded with $\bphi_F$, multiplied by the Green's function ($\Gv=\Lv^{-1}$), and then decoded with $\bpsi_u^{-1}$ to give an approximation of $\uv_k$. The $\Fv$ to $\uv$ cross-mapping loss is
    \begin{equation}
        \mathcal{L}_6 = \frac{1}{N} \sum_{k=1}^N \frac{\left\lVert \uv_k - \bpsi_u^{-1} \circ \Lv^{-1} \circ \bphi_F(\Fv_k) \right\rVert_2^2}{\left\lVert \uv_k\right\rVert_2^2}.
    \end{equation}
\end{enumerate}

Note that this final loss function gives the best indication of the performance of the network to solve the nonlinear BVP \eqref{BVP} using the Green's function. The strategy for solving \eqref{BVP} for a given discrete forcing function $\Fv$ is to encode the forcing function to obtain $\fv = \bphi_F(\Fv)$, apply the Green's function as in equation \eqref{Greens} to obtain $\vv$, and then decode this function to get the solution $\uv = \bpsi_u^{-1}(\vv)$. The discrete version of the convolution with the Green's function given in equation \eqref{Greens} is multiplication by the matrix $\Lv^{-1}$.
    
For the encoders $\bphi$ and $\bpsi$ and decoders $\bphi^{-1}$ and $\bpsi^{-1}$, we use a residual neural network (ResNet) architecture \cite{he2016deep}. The ResNet architecture has been successful in learning coordinate transformations for physical systems \cite{gin2020} and is motivated by near-identity transformations in physics. The linear operator $\Lv$ is constrained to be a real symmetric matrix and therefore is self-adjoint. Additionally, $\Lv$ is initialized as the identity matrix. Therefore, $\Lv$ is strictly diagonally dominant for at least the early parts of training which guarantees $\Lv$ is invertible and well-conditioned. For more information on the network architecture and training procedure, see Appendix \ref{sec:implementation}.

\section{Results}

The DeepGreen architecture, which is highlighted in Fig.~\ref{fig:greennet_arch} and whose detailed loss functions are discussed in the last section, is demonstrated on a number of canonical nonlinear BVPs. The first three BVPs are one-dimensional systems and the final one is a two-dimensional system. The nonlinearities in these problems do not allow for a fundamental solution, thus recourse is typically made to numerical computations to achieve a solution.  DeepGreen, however, can produce a fundamental solution which can then be used for any new forcing of the BVP. 

\subsection{Cubic Helmholtz} \label{sec:example}

The architecture and methodology is best illustrated using a basic example problem. The example problem uses a nonhomogeneous second-order nonlinear Sturm--Liouville model with constant coefficients and a cubic nonlinearity, thus making it a cubic Helmholtz equation.  The differential equation is given by
\begin{subequations}
\begin{gather} \label{eqn:model_1}
    u'' + \alpha u + \epsilon u^3 = F(x), \\
    u(0) = u(2\pi) = 0,
\end{gather}
\end{subequations}
where $u=u(x)$ is the solution when the system is forced with $F(x)$ with $x \in (0,2\pi)$, $\alpha = -1$ and $\epsilon = -0.3$. The notation $u''$ denotes $\frac{d^2}{dx^2}u(x)$. The dataset contains discretized solutions and forcings, $\{\uv_k,\Fv_k\}_{k=1}^N$. The forcing functions used for training are cosine and Gaussian functions; details of data generation and the forcing functions are provided in Appendix \ref{sec:datagen}.
\begin{figure}[t]
    \centering
    \includegraphics{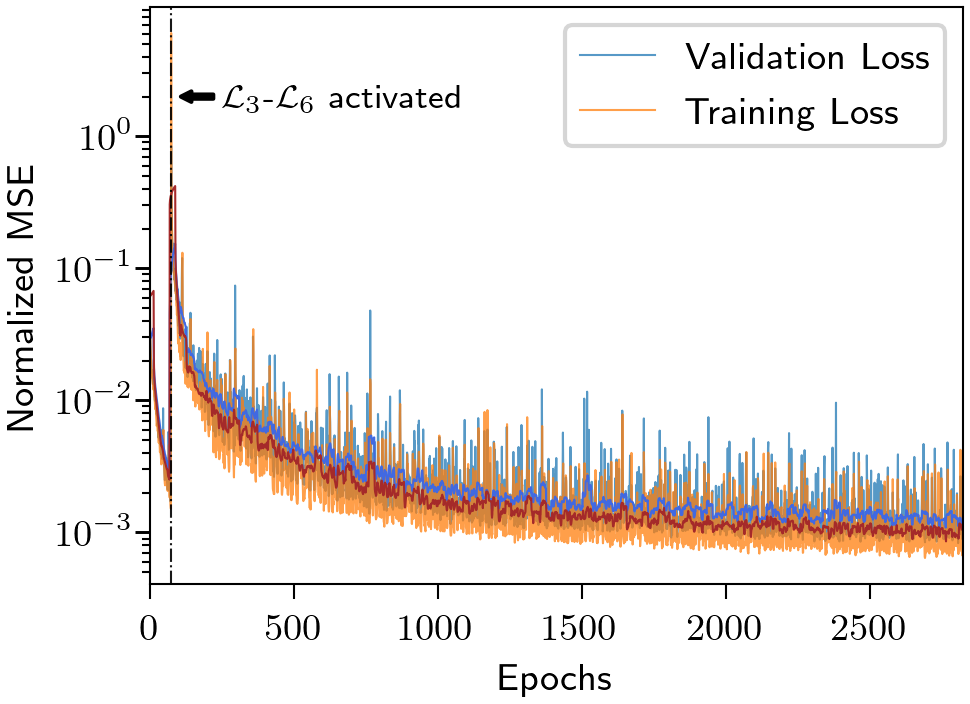}
    \caption{Learning curve. This is a typical learning curve for the DeepGreen  architecture. The vertical dashed line indicates where the training procedure transitions from autoencoders-only (only $\mathcal{L}_1$ and $\mathcal{L}_2$) to a full-network training procedure (all losses).}
    \label{fig:training_loss}
\end{figure}
\begin{figure}[t]
    \centering
    \includegraphics{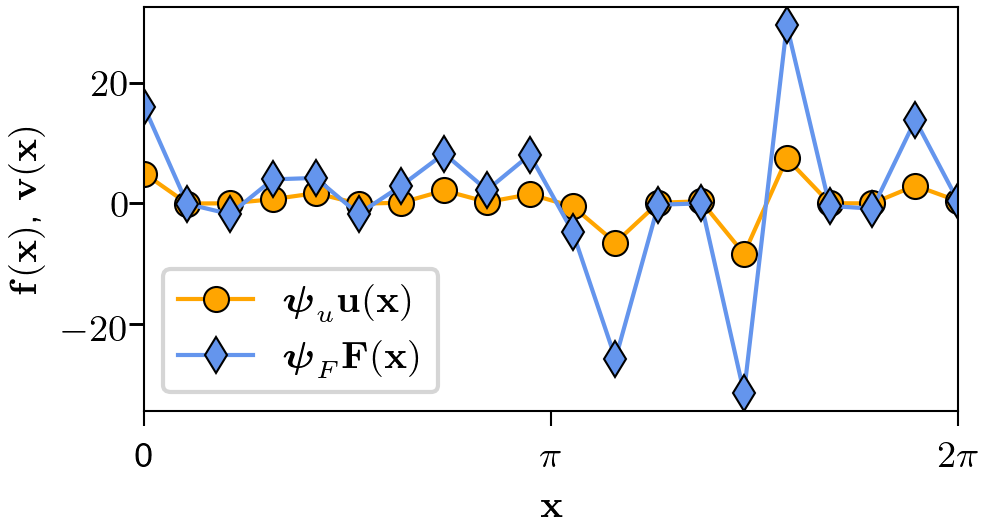}
    \caption{Latent space representations $\vv_k$ and $\fv_k$. The autoencoder transformation $\bpsi_u$ encodes $\uv_k$ to the latent space, producing the vector $\vv_k$ (orange). The forcing vector $\Fv_k$ is transformed by $\bpsi_F$ to the encoded vector $\fv_v$ (blue).}
    \label{fig:latent_rep}
\end{figure}
\begin{figure*}[t]
    \centering
    \includegraphics[width=12cm]{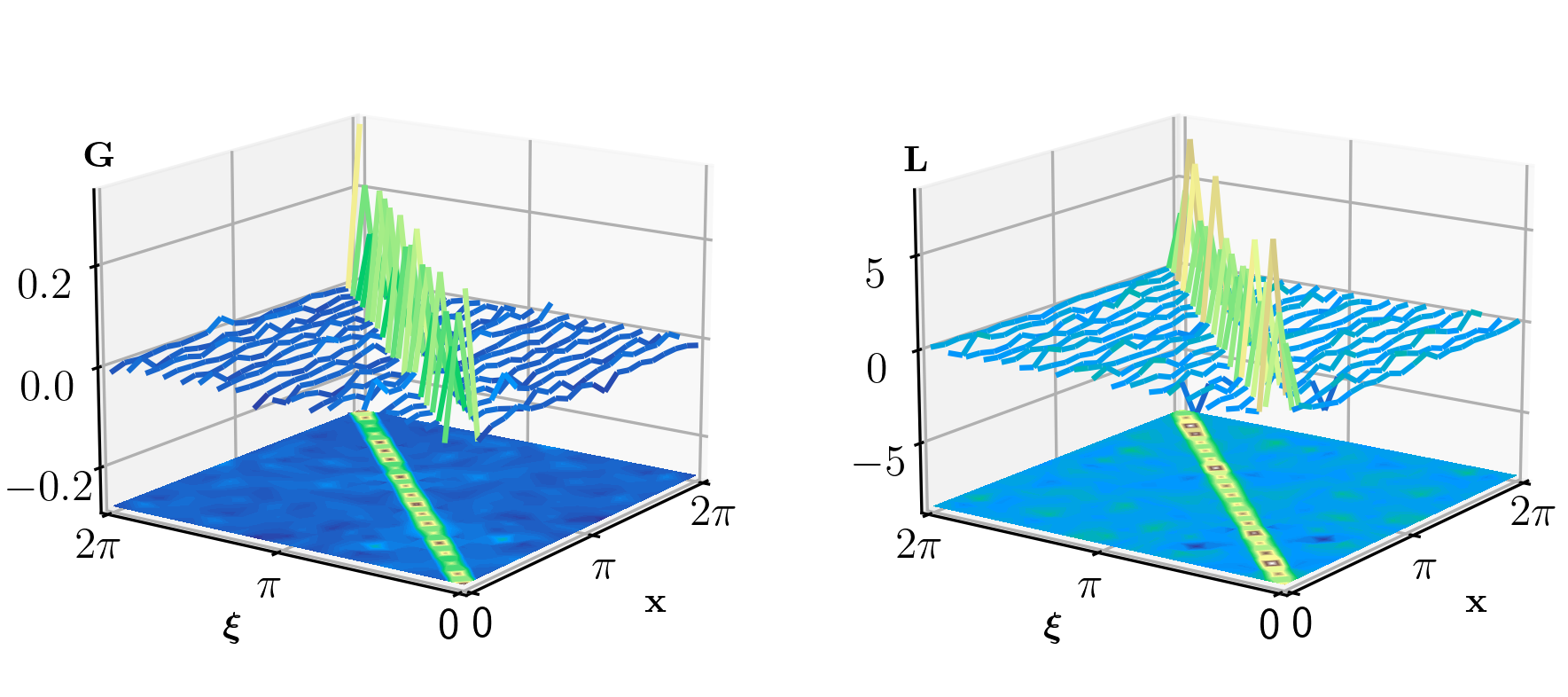}
    \vspace*{-.2in}
    \caption{Visualized operator and Green's function. Discovered Green's function $\Gv=\Lv^{-1}$ and corresponding linear operator $\Lv$.}
    \label{fig:LG_fig}
\end{figure*}
The data is divided into three groups: training, validation, and test. The training and validation sets are used for training the model. The test set is used to evaluate the results. The training set contains $N_{train}=8906$ vector pairs $\uv_k$ and $\Fv_k$. The validation set contains $N_{validation}=2227$ and test set contains $N_{test}=1238$.

\subsubsection{Training the Model}
The autoencoders used in this example are constructed with fully connected layers. In both autoencoders, a ResNet-like identity skip connection connects the input layer to the layer before dimension reduction in the encoder, and the first full-dimension layer in the decoder with the final output layer (see Figure \ref{fig:1d_layer_arch}).

The model is trained in a two-step procedure. First, the autoencoders are trained, without connection in the latent space, to condition the networks as autoencoders. In this first phase, only the autoencoder loss functions listed in Figure \ref{fig:greennet_arch} are active ($\mathcal{L}_1$ and $\mathcal{L}_2$). After a set number of epochs, the latent spaces are connected by an invertible matrix operator, $\Lv$, and the remaining 4 loss functions in Figure \ref{fig:greennet_arch} become active ($\mathcal{L}_3$--$\mathcal{L}_6$).
In the final phase of training, the autoencoder learns to encode a latent space representation of the system where properties  associated with linear systems hold true, such as linear superposition.

Figure \ref{fig:training_loss} shows a typical training loss curve. The vertical dashed line indicates the transition between the two training phases. The models in this work are trained for 75 epochs in the first autoencoder-only phase and 2750 epochs in the final phase. The first-phase epoch count was tuned empirically based on final model performance. The final phase epoch count was selected for practical reasons; the training curve tended to flatten around 2750 epochs in all of our tested systems.
The autoencoder latent spaces are critically important. The latent space is the transformed vector space where linear properties (e.g. superposition) are enforced which enables the solution of nonlinear problems. In the one-dimensional problems, the latent spaces vectors $\vv$ and $\fv$ are in $\mathbb{R}^{20}$.

The latent spaces did not have any obvious physical interpretation, and qualitatively appeared similar to the representations shown in Figure \ref{fig:latent_rep}. We trained 100 models to check the consistency in the learned model and latent space representations, but discovered the latent spaces varied considerably (see Appendix \ref{sec:additional_results}). This implies the existence of an infinity of solutions to the coordinate transform problem, which indicates further constraints could be placed on the model.

\begin{figure*}[t]
    \centering
    \includegraphics[width=13cm]{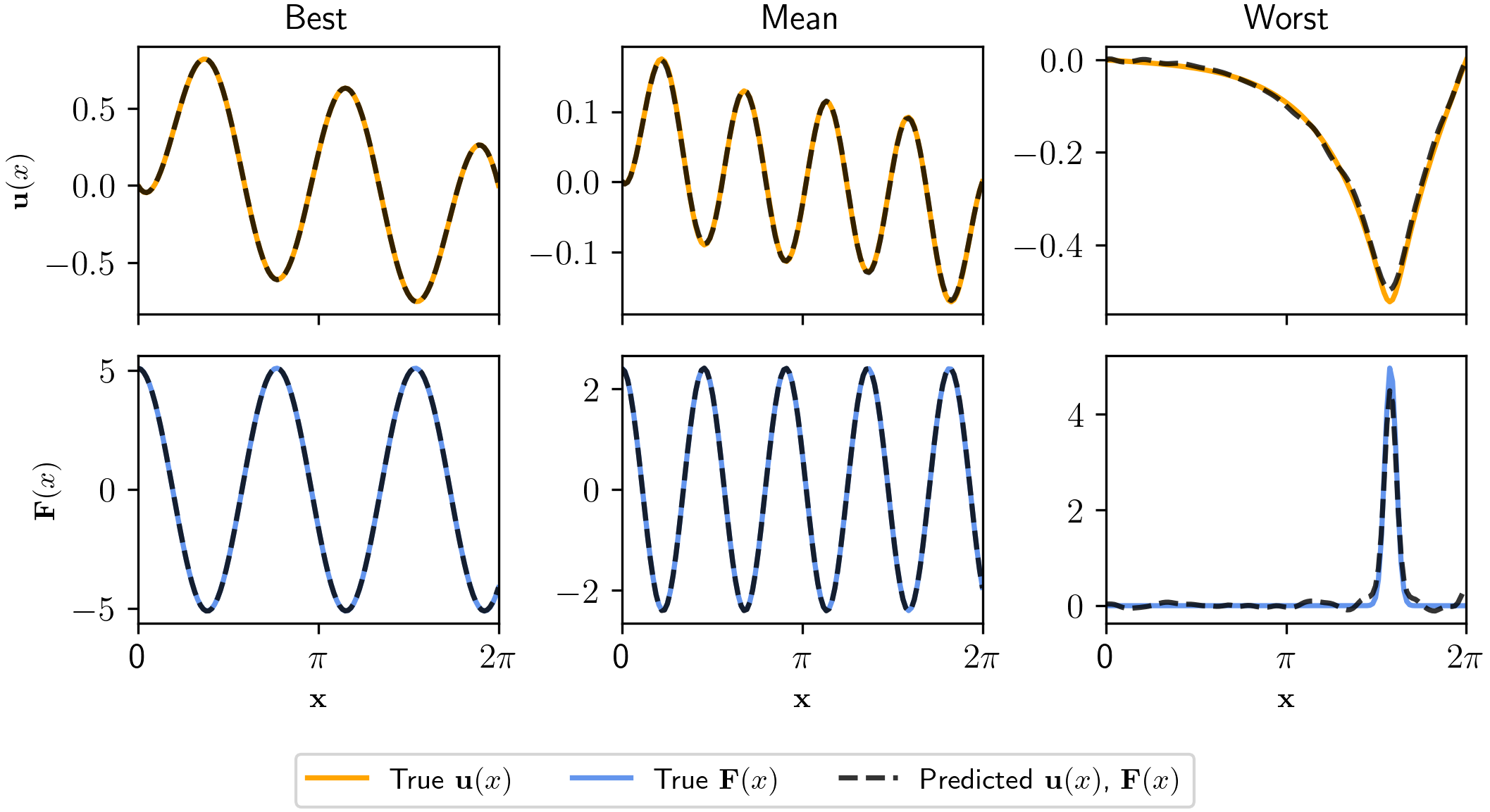}
    \caption{Model predictions on test data. The top row shows the true solution $\uv_k(x)$ and the solution predicted by the network given the forcing $\Fv_k(x)$ using the Green's function $\Gv$. The bottom row shows the true forcing function $\Fv_k(x)$ compared to the forcing computed by applying the operator $\Lv$ to the solution $\uv_k$. Three columns show the best, mean, and worst case samples as evaluated by the sum of normalized $\ell 2$ reconstruction errors.}
    \label{fig:similar_test}
\end{figure*}

Despite lacking obvious physical interpretations, the latent space enables discovery of an invertible operator $\Lv$ which described the linear system $\Lv[\vv_k]=\fv_k$. The operator matrix $\Lv$ can be inverted to yield the Green's function matrix $\Gv$, which allows computation of solutions to the linearized system $\vv_k = \Gv[\fv_k]$. An example of the operator $\Lv$ and its inverse $\Gv$ are shown in Figure \ref{fig:LG_fig}.
The operator and Green's function shown in Figure \ref{fig:LG_fig} display an important prominent feature seen in all of the results: a diagonally-dominant structure. We initialize the operator as an identity matrix, but the initialization had little impact on the diagonally-dominant form of the learned operator and Green's function matrices (see Appendix \ref{sec:additional_results}). The diagonally-dominant operators indicate that the deep learning network tends to discover a coordinate transform yielding a nearly-orthonormal basis, which mirrors the common approach of diagonalization in spectral theory for Hermitian operators. Furthermore, diagonally-dominant matrices guarantee favorable properties for this application such as being well-conditioned and non-singular.

We emphasize that training parameters and model construction choices used in this work were not extensively optimized. We expect the model performance can be improved in a myriad of ways including extending training times, optimizing model architecture, modifying the size of the latent spaces, restricting the form of the operator, and applying additional constraints to the model. However, these topics are not the main scope of the present work; our focus is to illustrate the use of autoencoders as a coordinate transform for finding solutions to nonlinear BVPs.
\subsubsection{Evaluating the Model}
The goal for this model is to find a Green's function $\Gv$ for computing solutions $\uv_k$ to a nonlinear BVP governed by  (\ref{BVP_continuous}) for a given forcing function $\Fv_k$. Similarly, we can estimate the forcing term, $\Fv_k$, given the solution $\uv_k$. The model is consequently evaluated by its ability to use the learned Green's function and operator for predicting solutions and forcings, respectively, for new problems from a withheld test data set.

Recall the original model is trained on data where the forcing function is a cosine or Gaussian function.  As shown in Figure \ref{fig:similar_test}, the model performs well on withheld test data where the forcing functions are cosine or Gaussian functions, producing a cumulative loss around $10^{-4}$. The solutions $\uv_k$ and forcing $\Fv_k$ are depicted for the best, mean, and worst samples scored by cumulative loss.

It's important to note the test data used in Figure \ref{fig:similar_test} is similar to the training and validation data. Because ML models typically work extremely well in interpolation problems, it is reasonable to expect the model to perform well on this test data set.

\begin{figure*}[t]
    \centering
    \includegraphics[width=13cm]{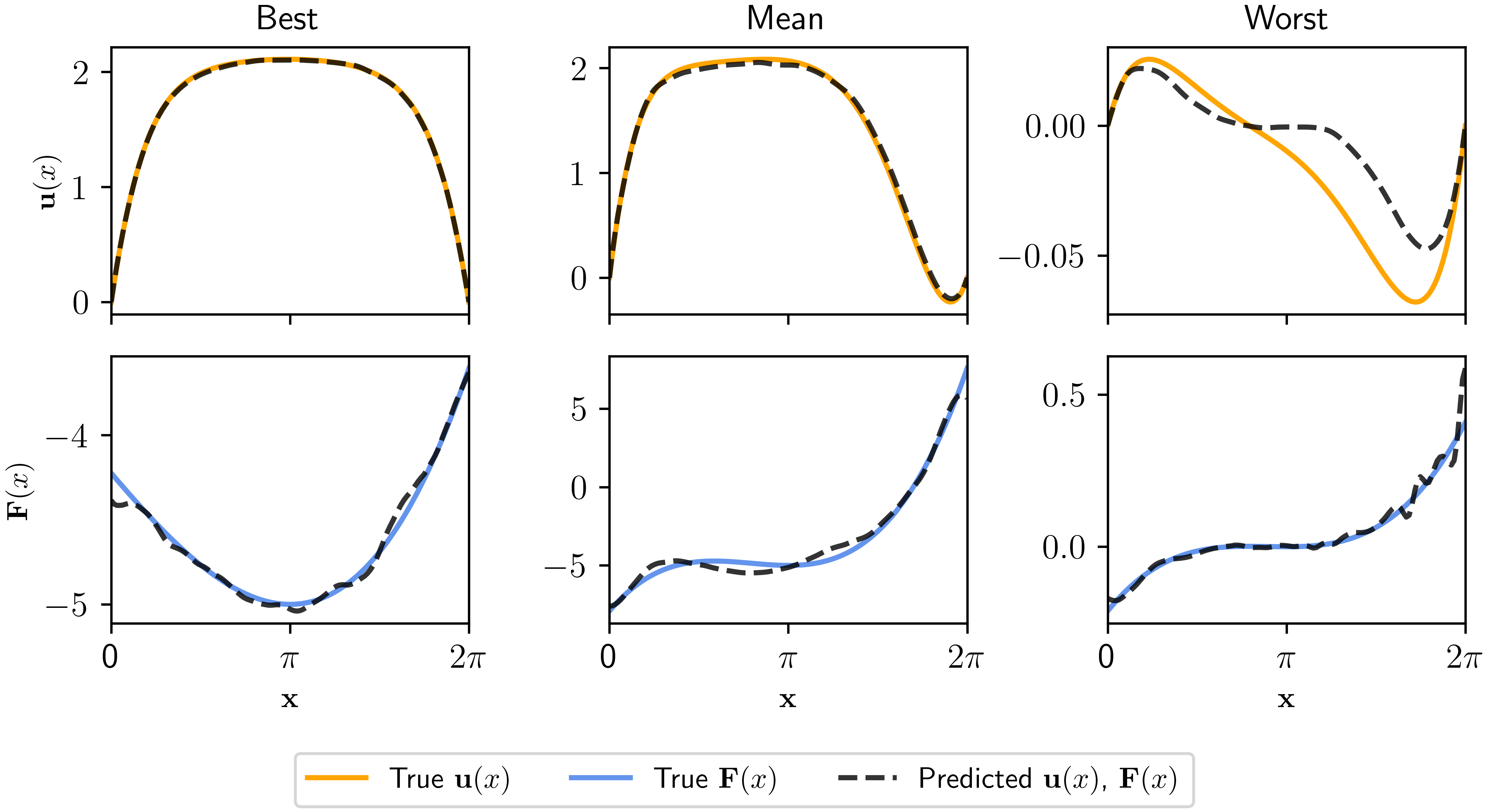}
    \caption{Model predictions on cubic Helmholtz forced system. The top row shows the true solution $\uv_k(x)$ and the solution predicted by the network given the forcing $\Fv_k(x)$ using the Green's function $\Gv$. The bottom row shows the true forcing function $\Fv_k(x)$ compared to the forcing computed by applying the operator $\Lv$ to the solution $\uv_k$. Three columns show the best, mean, and worst case samples as evaluated by the sum of normalized $\ell 2$ reconstruction errors.}
    \label{fig:dissimilar_test}
\end{figure*}

As an interesting test to demonstrate the ability of the model to extrapolate, we prepared a separate set of test data $\{\uv_k,\Fv_k\}_{k=1}^N$ containing solutions where $\Fv_k$ are cubic polynomial forcing functions. This type of data was not present in training, and provides some insight into the generality of the learned linear operator and Green's function matrices. Figure \ref{fig:dissimilar_test} shows examples of how the model performs on these cubic polynomial-type forcing functions. Similar to Figure \ref{fig:similar_test}, the best, mean, and worst samples are shown as graded by overall loss.
\begin{figure*}[t]
    \centering
    \includegraphics[width=13cm]{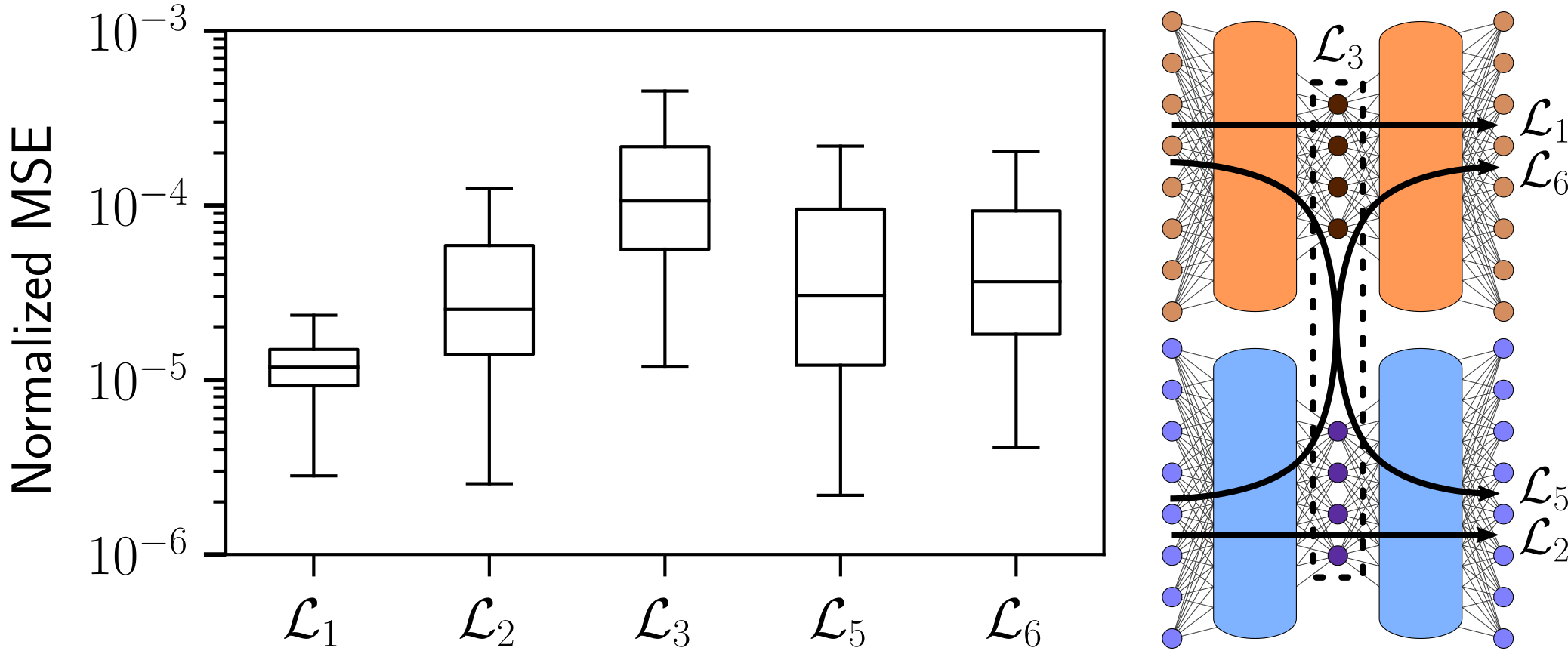}
    \caption{Model performance summary. Distribution of loss values are shown for every sample in the test data set. Model loss functions are minimized during training, making them a natural metric to use for summarizing performance.}
    \label{fig:box plot}
\end{figure*}
Figures \ref{fig:similar_test} and \ref{fig:dissimilar_test} provide some qualitative insight into the model's performance on specific instances selected from the pool of evaluated data. A quantitative perspective of the model's performance is presented in Figure \ref{fig:box plot}. This box plot shows statistics (median value, $Q_1$, $Q_3$, and range) for four of the loss functions evaluated on the similar (cosine and Gaussian) test data. Note the superposition loss function is \emph{not} scored in this plot because the superposition loss function can only be evaluated within a single batch, and the loss depends on batch size and composition.

In conclusion, the DeepGreen architecture enables  discovery of invertible, linearizing transformations that facilitate identification of a linear operator and Green's function to solve nonlinear BVPs. It is tested on data similar and dissimilar to the training data, and evaluated on the loss functions that guide the training procedure. The discovered operator and Green's function take on a surprisingly diagonally-dominant structure, which hints at the model's preference to learn an optimal basis. The model appears to extrapolate beyond the test data, suggesting that the learned operator is somewhat general to the system.

\subsection{Nonlinear Sturm--Liouville and Biharmonic Operators}
In addition to the example system described above, the approach was applied to two other one-dimensional systems. We used the same training procedure and forcing functions that were described in Section \ref{sec:example}. The first is a system governed by the nonlinear Sturm--Liouville equation 
\begin{gather*}
    [-p(x) u']' + q(x) (u + \epsilon u^3) = F(x), \\
    u(0) = u(2\pi) = 0,
\end{gather*}
where $\epsilon = 0.4$ controls the extent of nonlinearity, and $p(x)$ and $q(x)$ are spatially-varying coefficients
\begin{gather*}
    p(x) = 0.5 \sin(x) - 3, \\
    q(x) = 0.6 \sin(x) - 2 ,
\end{gather*}
with $x \in [0,2\pi]$. The final one-dimensional system is a biharmonic operator with an added cubic nonlinearity
\begin{gather*}
    [-p u'']'' + q (u + \epsilon u^3) = F(x), \\
    u(0) = u(2\pi) = u'(0) = u'(2\pi) = 0,
\end{gather*}
where $p=-4$ and $q=2$ are the coefficients and $\epsilon=0.4$ controls the nonlinearity. As in the prior example, the forcing functions in the training data are cosine and Gaussian functions, which are described further in Appendix \ref{sec:datagen}.

Results for all the one-dimensional models, including the cubic Helmholtz example from Section \ref{sec:example}, are presented in Table \ref{tab:summary_table}. Model performance is quantitatively summarized by box plots and the Green's function matrix is shown for each model.

\begin{table*}[ht!]
    \centering
    {\renewcommand{\arraystretch}{1.3}
    \begin{tabular}
    	{  | c | c | c | }
        \hline

    	\makecell{Nonlinear\\ 
    	cubic Helmholtz \\ 
    	\small (constant coefficients) \\ 
    	$u'' + \alpha u + \epsilon u^3 = F$ \\ 
    	$u(0) = u(L) = 0$} &
    	
    	\makecell{Nonlinear\\
    	Sturm--Liouville \\ 
    	\small (varying $p(x)$, $q(x)$) \\ 
    	$[-p u']' \! + \! q u \! + \! \alpha q u^3 \! = \! F$ \\
    	$u(0) = u(L) = 0$} & 
    	
    	\makecell{Nonlinear \\ 
    	Biharmonic operator \\
    	\small (constant coefficients) \\ 
    	$-p u'''' \! + \! q u \! + \! \alpha q u^3 \! = \! F$ \\ 
    	$u(0) = u(L) = 0$} 
    	
    	 \\ \hline
    	 
    	 \makecell{\centering \includegraphics[width=47mm]{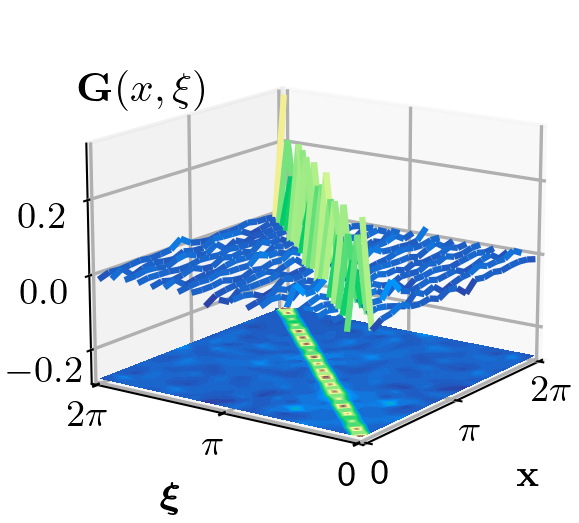}} & 
    	 
    	 \makecell{\centering \includegraphics[width=47mm]{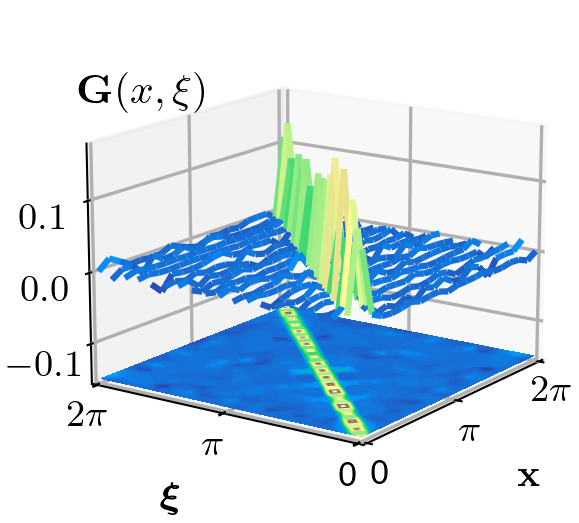}} & 
    	 
    	 \makecell{\centering \includegraphics[width=47mm]{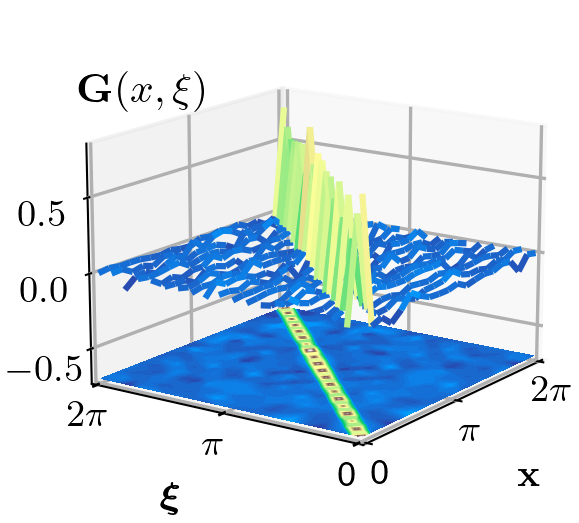}}
    	 
    	 \\ \hline

    	 \multicolumn{3}{| c |}{\makecell{\\
    	 \includegraphics[width=130mm]{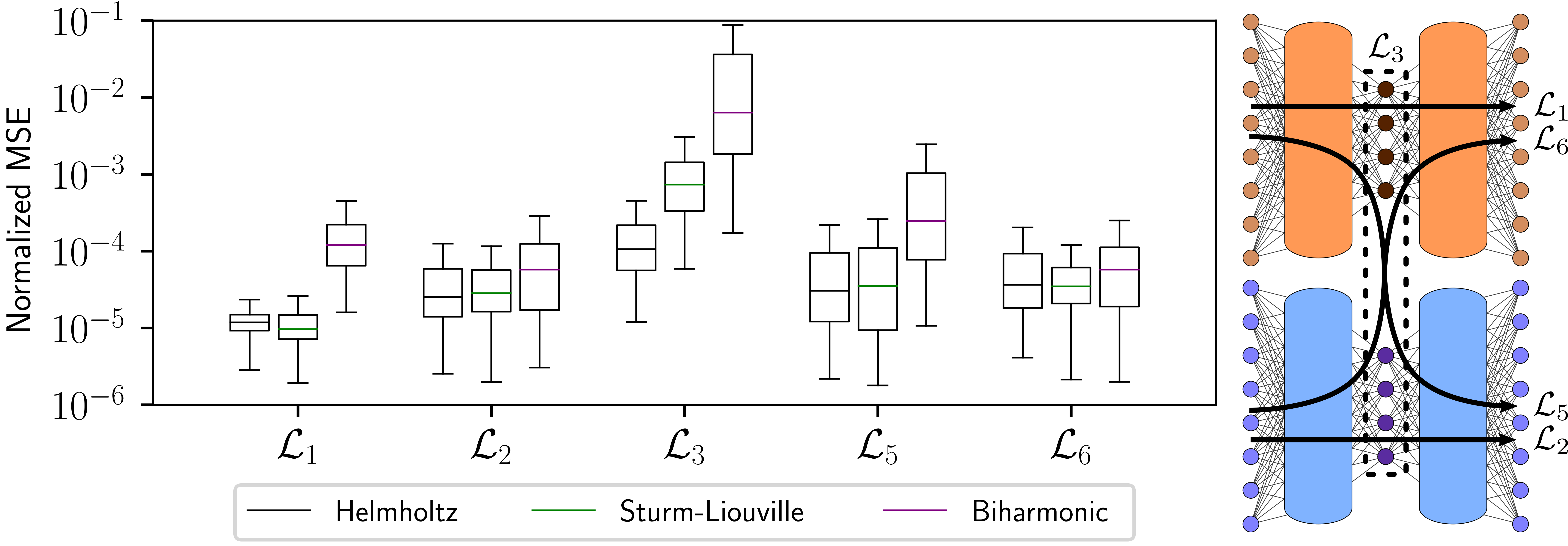}}}
    	 
    	 \\ \hline

    \end{tabular}
    }
    \caption{Summary of results for three one-dimensional models. The models are provided with the Green's function learned by DeepGreen. A summary box plot shows the relative losses $\mathcal{L}_1$, $\mathcal{L}_2$, $\mathcal{L}_3$, $\mathcal{L}_5$, and $\mathcal{L}_6$ for all three model systems.}
    \label{tab:summary_table}
\end{table*}

Importantly, the learned Green's function matrices consistently exhibit diagonally-dominant structure. The losses for the nonlinear cubic Helmholtz equation and the nonlinear Sturm--Liouville equation are similar which indicates that spatially-varying coefficients do not make the problem significantly more difficult for the DeepGreen architecture. In contrast, the loss for the nonlinear biharmonic equation are about an order of magnitude higher than the other two systems. This result implies the fourth-order problem is more difficult than the second-order problems. Also of note is that the linear operator loss $\mathcal{L}_3$ is consistently the highest loss across all models. Therefore, it is easier for DeepGreen to find invertible transformations for the solutions and forcing functions than it is to find a linear operator that connects the two latent spaces.

\subsection{Nonlinear Poisson Equation}

\begin{figure*}[t]
\vspace*{-.15in}
\centering
\includegraphics[width=15cm]{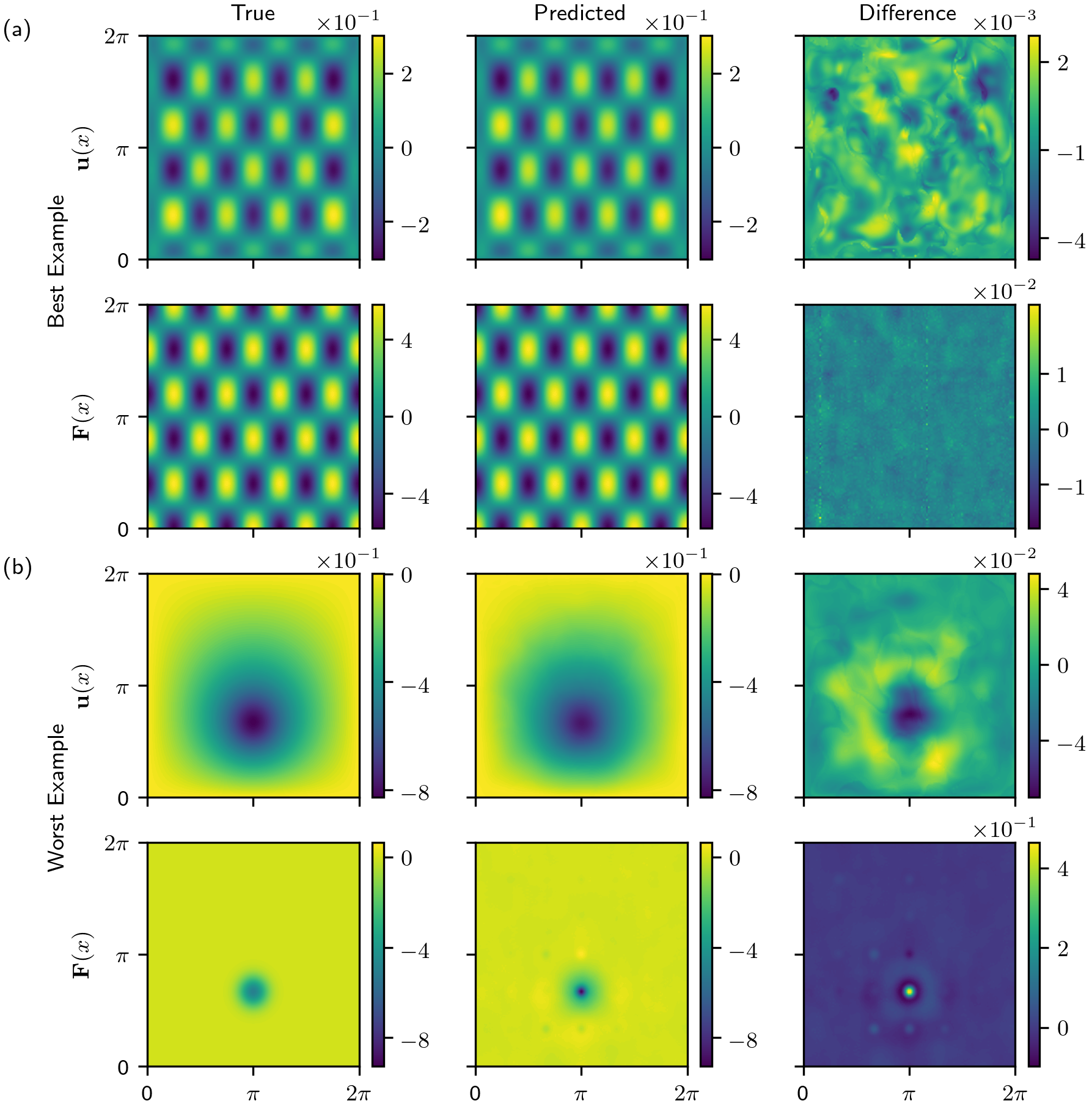}
\vspace*{-.15in}
\caption{Model predictions for the (a) best and (b) worst examples from test data with Gaussian and cosine forcings. In both (a) and (b), the top row shows the true solution $\uv(\bx)$, the predicted solution using the Green's function, and the difference between the true and predicted solution. The bottom row shows the true forcing function $\Fv(\bx)$, the predicted forcing function, and the difference between the true and predicted forces. In order to account for the difference in scale between $\uv(\bx)$ and $\Fv(\bx)$, the differences are scaled by the infinity norm of the true solution or forcing function ($\text{Difference} = (\text{True} - \text{Predicted})/ || \text{True}||_{\infty}$).}
\label{fig:similar_test_2D}
\end{figure*}

We also tested our method on a two-dimensional system. The two-dimensional model is a nonlinear version of the Poisson equation with Dirichlet boundary conditions
\begin{subequations} \label{NLP}
\begin{align}
    &-\nabla \cdot \left[(1+u^2) \nabla u\right] = F(\xv), &
    \bx &\in \Omega,  \\
     &u = 0, &
     \bx &\in \partial\Omega,
\end{align}
\end{subequations}
where $\Omega := (0,2\pi) \times (0,2 \pi)$. Similar to the one-dimensional models, the forcing functions used to train the model are cosine and Gaussian functions, the details of which are provided in Appendix \ref{sec:datagen}. The sizes of the data sets are also similar to the one-dimensional data sets. The training data contains $N_{train}=9806$ vector pairs $\uv_k$ and $\Fv_k$, the validation data contains $N_{validation}=2452$, and the test data contains $N_{test}=1363$. 

The network architecture of the encoders and decoders for the two-dimensional example differs from the one-dimensional examples. Instead of fully connected layers, convolutional layers were used in the encoders and decoders. However, we still use a ResNet architecture. Additionally, the latent space vectors are in $\mathbb{R}^{200}$. Full details on the network architecture can be found in Appendix \ref{sec:implementation}. Note that the method proposed for discovering Green's functions allows for any network architecture to be used for the encoders and decoders. For the one-dimensional example, similar results were obtained using fully connected and convolutional layers. However, the convolutional architecture was better in the two-layer case and also allowed for a more manageable number of parameters for the wider network that resulted from discretizing the two-dimensional space.

\begin{figure*}[t]
\centering
\includegraphics[width=15cm]{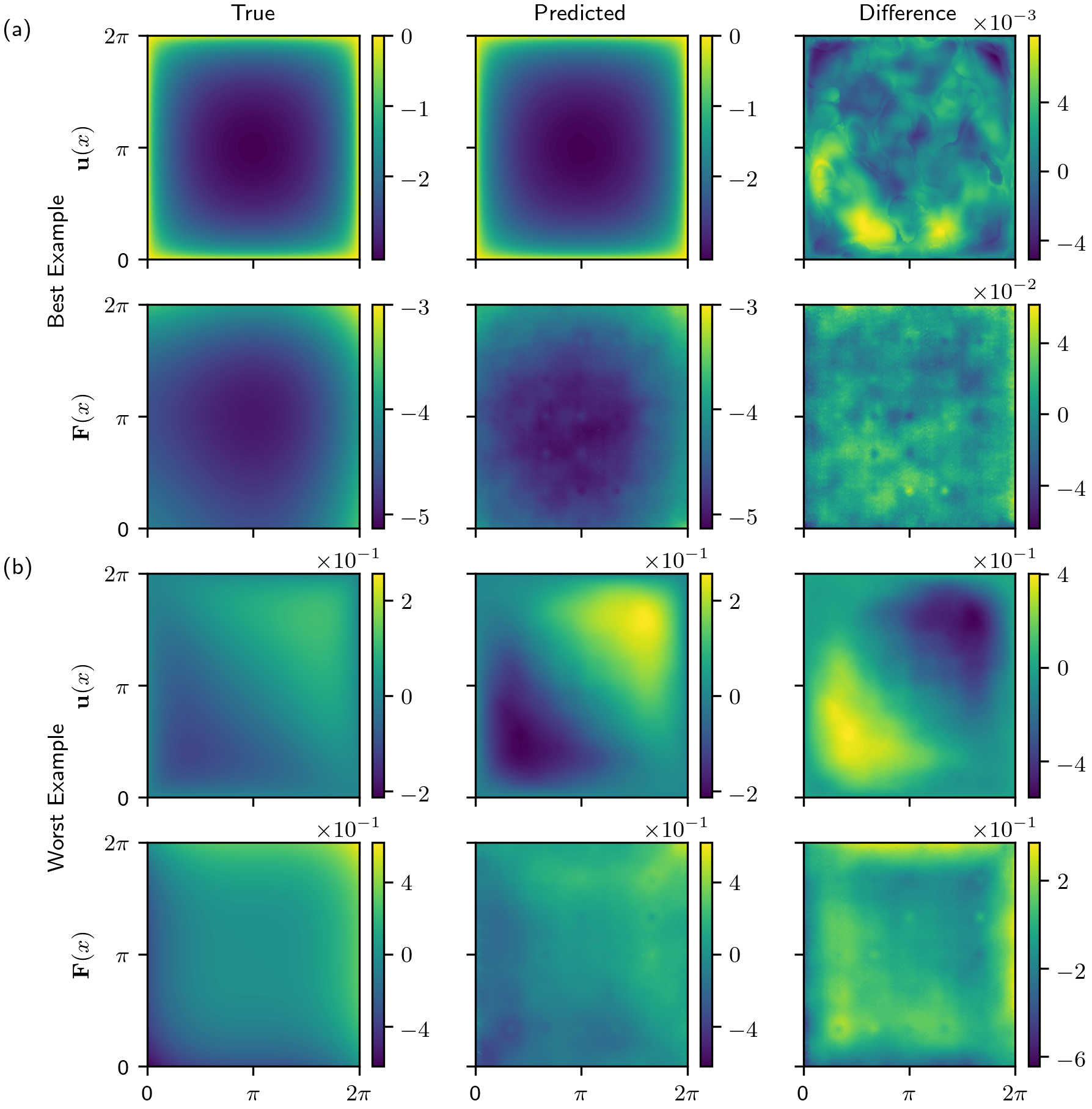}
\vspace{-.15in}
\caption{Model predictions for the (a) best and (b) worst examples from test data with cubic polynomial forcings. In both (a) and (b), the top row shows the true solution $\uv(\bx)$, the predicted solution using the Green's function, and the difference between the true and predicted solution. The bottom row shows the true forcing function $\Fv(\bx)$, the predicted forcing function, and the difference between the true and predicted forces. In order to account for the difference in scale between $\uv(\bx)$ and $\Fv(\bx)$, the differences are scaled by the infinity norm of the true solution or forcing function ($\text{Difference} = (\text{True} - \text{Predicted})/ || \text{True}||_{\infty}$).}
\label{fig:dissimilar_test_2D}
\end{figure*}

The operator and Green's function for the two-dimensional model are similar to those displayed in shown in Figure \ref{fig:LG_fig}.
The diagonal dominance is even more prevalent in this case than the one-dimensional example.
The model was evaluated on test data containing cosine and Gaussian forcing functions. Figure \ref{fig:similar_test_2D}a shows the true solution $\uv(x)$ and forcing function $\Fv(x)$ as well as the network predictions for the example from the test data for which the model performed the best (i.e. the smallest value of the loss). The difference between the true and predicted functions is shown in the right column of Figure \ref{fig:similar_test_2D}a and is scaled by the infinity norm of the true solution or forcing functions. Figure \ref{fig:similar_test_2D}b shows similar results but for the worst example from the test data. In both cases, the model gives a qualitatively correct solution for both $\uv(x)$ and $\Fv(x)$. Unsurprisingly, the network struggles most on highly localized forcing functions and has the highest error in the region where the forcing occurs.

\begin{figure*}[t]
    \centering
    \includegraphics[width=13cm]{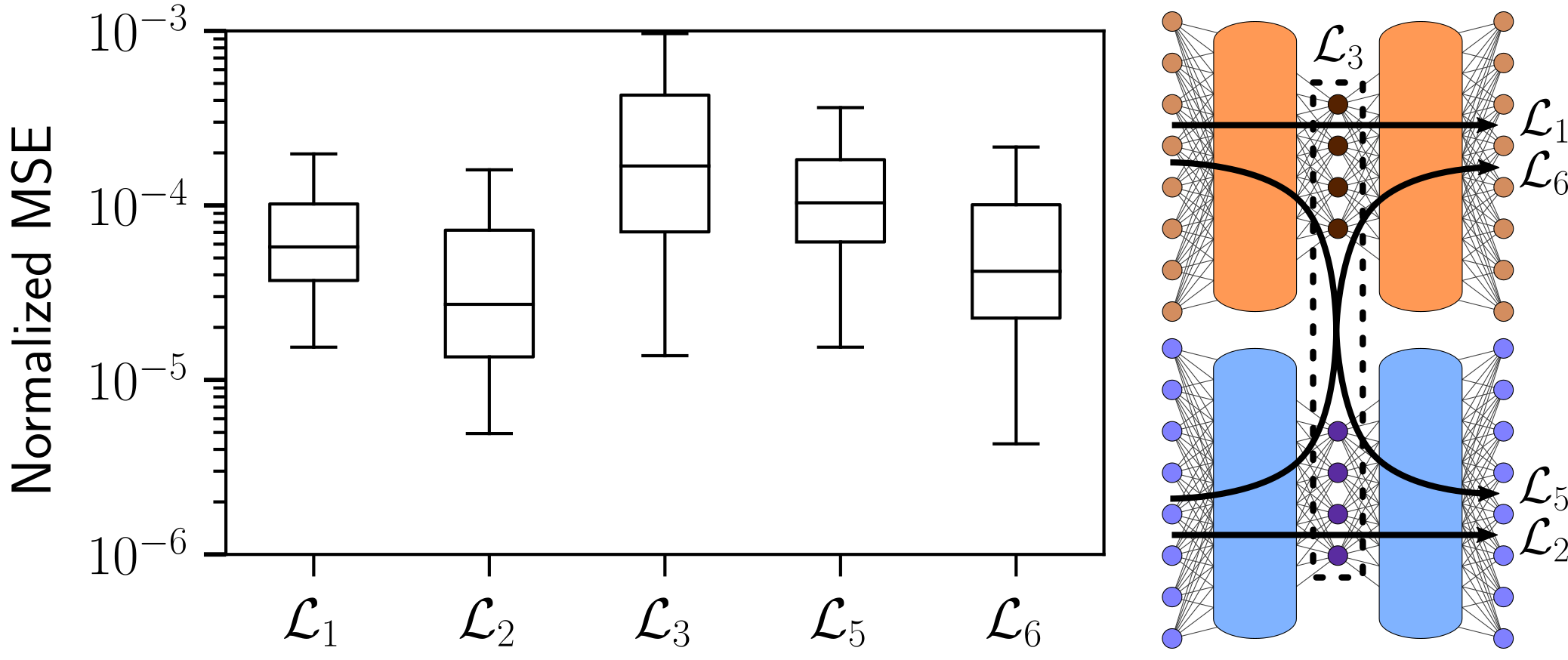}
    \caption{Two-dimensional Poisson model performance summary. Distribution of loss values are shown for every sample in the test data set. Model loss functions are minimized during training, making them a natural metric to use for summarizing performance.}
    \label{fig:boxplot_2D}
\end{figure*}

The model was also evaluated on test data that has cubic polynomial forcing functions, a type of forcing function not found in the training data. The best and worst examples are shown in Figure \ref{fig:dissimilar_test_2D}. Although the model does not perform as well for test data which is not similar to the training data, the qualitative features of the predicted solutions are still consistent with the true solutions.
Figure \ref{fig:boxplot_2D} shows a box plot of the model's performance on the similar (cosine and Gaussian forcing test data). The results are similar to the one-dimensional results, and, in fact, better than the biharmonic operator model.

\section{Conclusion}

We have leveraged the expressive capabilities of deep learning to discover linearizing coordinates for nonlinear BVPs, thus allowing for the construction of the {\em fundamental solution or nonlinear Green's function}.  Much like the Koopman operator for time-dependent problems, the linearizing transformation provides a framework whereby the fundamental solution of the linear operator can be constructed and used for any arbitrary forcing.  This provides a broadly applicable mathematical architecture for constructing solutions for nonlinear BVPs, which typically rely on numerical methods to achieve solutions.  Our DeepGreen architecture can achieve solutions for arbitrary forcings by simply computing the convolution of the forcing with the Green's function in the linearized coordinates.  

Given the critical role that BVPs play in the mathematical analysis of constrained physical systems subjected to external forces, the DeepGreen architecture can be broadly applied in nearly every engineering discipline since BVPs are prevalent in diverse problem domains including fluid mechanics, electromagnetics, quantum mechanics, and elasticity.  Importantly, DeepGreen provides a bridge between a classic and widely used solution technique to nonlinear BVP problems which generically do not have principled techniques for achieving solutions aside from brute-force computation.  DeepGreen establishes this bridge by providing a transformation which allows linear superposition to hold.  DeepGreen is a flexible, data-driven, deep learning approach to solving nonlinear boundary value problems (BVPs) using a dual-autoencoder architecture.
The autoencoders discover an invertible coordinate transform that linearizes the nonlinear BVP and identifies both a linear operator $L$ and Green's function $G$ which can be used to solve new  nonlinear BVPs.
We demonstrated that the method succeeds on a variety of nonlinear systems including nonlinear Helmholtz and Sturm--Liouville problems, nonlinear elasticity, and a 2D nonlinear Poisson equation.
The method merges the strengths of the universal approximation capabilities of deep learning with the physics knowledge of Green's functions to yield a flexible tool for identifying fundamental solutions to a variety of nonlinear systems.

\section*{Acknowledgments}
SLB is grateful for funding support from the Army Research Office (ARO W911NF-17-1-0306). JNK acknowledges support from the Air Force Office of Scientific Research (FA9550-19-1-0011).

\begin{spacing}{.01}
	\small
    \bibliographystyle{ScienceAdvances}
    \bibliography{references}
\end{spacing}

\newpage

\begin{appendices}

\section{Data Generation} \label{sec:datagen}

\subsection{1D Problems}
The data for all of the one-dimensional systems are created using the same method and forcing functions. Each solution is computed on an evenly-spaced 128-point grid using MATLAB's bvp5c solver with a relative error tolerance of $10^{-8}$ and an absolute error tolerance of $10^{-10}$. The forcing functions $\Fv_k(x)$ are designed to yield a variety of solutions $\uv_k$ such that $||\uv_k||_2 \simeq 1$. 

The training data consists of two types of systems: Gaussian-forced and cosine-forced systems. The Gaussian-forced systems have forcing functions of the form
\begin{equation*}
    F_k(x) = a \exp\left(\frac{- (x-b)^2}{2c^2}\right),
\end{equation*}
where $a \in \{-25, -20, -15, -10, -5, 5, 10, 15, 20, 25\}$, $b \in \{0, 2\pi/19, 4\pi/19, \dots, 2\pi\}$, and $c \in \{0.1, 0.3, 0.5, \dots, 4.9\}$. The cosine forcing functions are of the form
\begin{equation*}
    F_k(x) = \alpha \cos(\beta x),
\end{equation*}
where $\alpha \in \{1, 1.1, 1.2, \dots, 10\}$ and $\beta \in \{1, 1.05, 1.10, \dots, 5\}$. This gives a total of 5000 Gaussian-forced solutions and 7371 cosine-forced solutions. For the cubic Helmholtz equation and the nonlinear Sturm--Liouville equation with spatially-varying coefficients, all of the 12371 solutions were within the error tolerance. However, there were 97 solutions of the nonlinear biharmonic equation that did not meet the error tolerance and were therefore discarded. Of the remaining data, 10\% are randomly chosen and withheld as test data, 80\% are used as training data, and 20\% are used as validation data.

In order to test the ability of the network to generalize, we also have another test data set that consists of solutions with
cubic forcing functions of the form
\begin{equation*}
  F_i(x) = \gamma (x-\pi)^3,
\end{equation*}
where $\gamma \in \{0.01, 0.03, 0.05, \dots, 0.29\}$, and cubic forcing functions of the form
\begin{equation*}
  F_i(x) = \gamma (x-\pi)^3 + \zeta (x-\pi)^2 + \bpsi,
\end{equation*}
where $\gamma \in \{0.01, 0.03, 0.05, \dots, 0.29\}$, $\zeta \in   \{0.01, 0.03, 0.05, \dots, 0.49\}$, and $\bpsi \in \{-5, -4, -3, \dots, 5\}$. There are a total of 4140 solutions with cubic forcing functions.

\subsection{2D Problem}
The two-dimensional data satisfies the nonlinear Poisson equation \eqref{NLP}. The solutions are computed with a finite element method using the DOLFIN library \cite{LoggWells2010a} of the FEniCS Project \cite{AlnaesBlechta2015a, LoggMardalEtAl2012a}. The forcing functions are similar to the one-dimensional data in that there are Gaussian and cosine forcing functions along with a separate data set of cubic polynomial forcing functions used to test the ability of the network to generalize. The Gaussian forcing functions are of the form
\begin{equation*}
    F_k(x,y) = a \exp\left(\frac{- (x-b_x)^2 - (y-b_y)^2}{2c^2}\right),
\end{equation*}
where $a \in \{-25, -20, -15, -10, -5, 5, 10, 15, 20, 25\}$, $b_x, b_y \in \{\pi/3, 2\pi/3, \pi, 4\pi/3, 5\pi/3\}$, and $c \in \{0.1, 0.3, 0.5, \dots, 4.9\}$. The cosine forcing functions are of the form 
\begin{equation*}
    F_k(x,y) = \alpha \cos(\beta_x x) \cos(\beta_y y),
\end{equation*}
where $\alpha \in \{1, 1.1, 1.2, \dots, 10\}$ and $\beta_x, \beta_y \in \{1, 1.5, 2, \dots, 5\}$. The cubic forcing functions are of the form
\begin{equation*}
  F_i(x, y) = \gamma_x (x-\pi)^3 + \gamma_y (y-\pi)^3,
\end{equation*}
where $\gamma_x, \gamma_y \in \{0.01 + 0.28k/3 | k = 0,1,2,3\}$, and cubic forcing functions of the form
\begin{equation*}
\begin{split}
  F_i(x,y) = &\gamma_x (x-\pi)^3 + \gamma_y (y-\pi)^3 \\ &+ \zeta_x (x-\pi)^2 +  \zeta_y (y-\pi)^2 + \bpsi,
\end{split}
\end{equation*}
where $\gamma_x, \gamma_y \in \{0.01 + 0.28k/3 | k = 0,1,2,3\}$, $\zeta_x, \zeta_y \in   \{0.01, 0.07, 0.13, 0.19, 0.25\}$, and $\bpsi \in \{-5, -4, -3, \dots, 5\}$. There are 6250 solutions with Gaussian forcing functions, 7371 solutions with cosine forcing functions, and 4416 solutions with cubic forcing functions.

\begin{figure*}[t]
    \centering
    \includegraphics{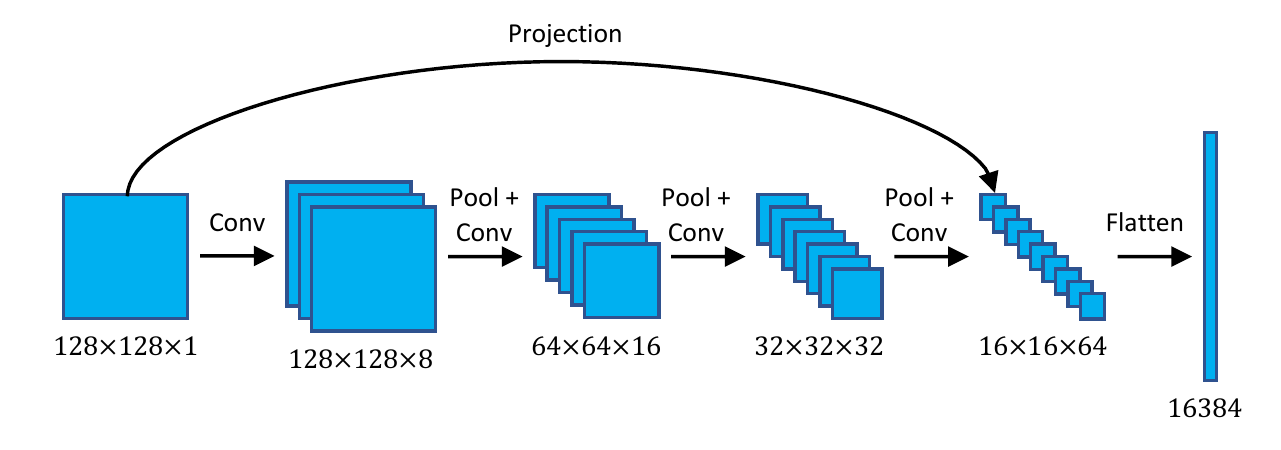}
    \vspace{-.2in}
    \caption{Encoder network architecture for the two-dimensional data. All convolutional layers use $4 \times 4$ kernels with stride size 1, zero-padding, and ReLU activation functions. All pooling layers are average pooling layers with pool size 2 and stride size 2.}
    \label{fig:2D_architecture}
\end{figure*}

\begin{figure*}[t]
    \centering
    \includegraphics{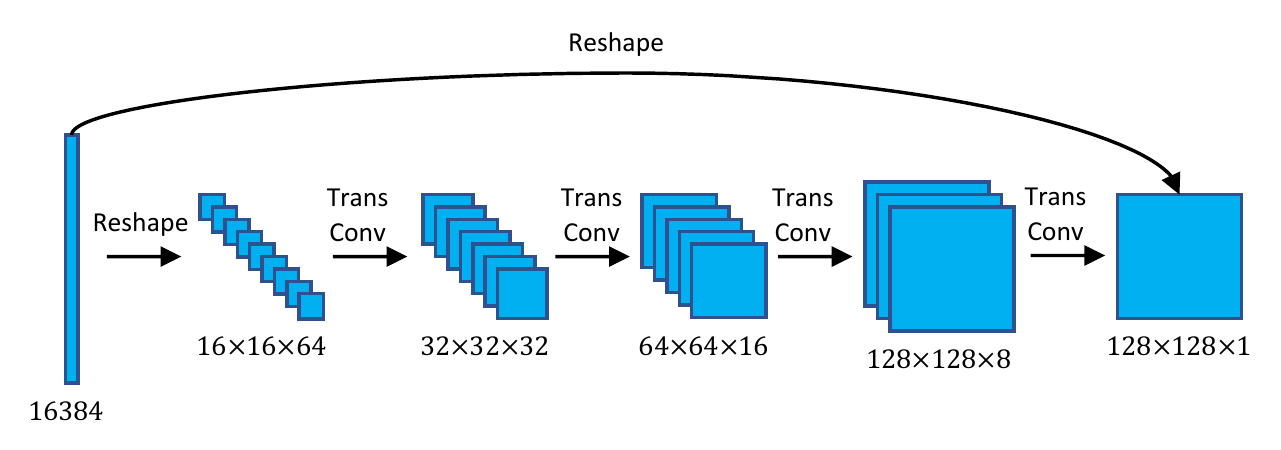}
    \vspace{-.2in}
    \caption{Decoder network architecture for the two-dimensional data. All transposed convolutional layers use $4 \times 4$ kernels with stride size 2, zero-padding, and ReLU activation functions except for the last layer which has stride size 1.}
    \label{fig:2D_architecture2}
\end{figure*}

\begin{figure}[t]
    \centering
    \includegraphics[width=25mm]{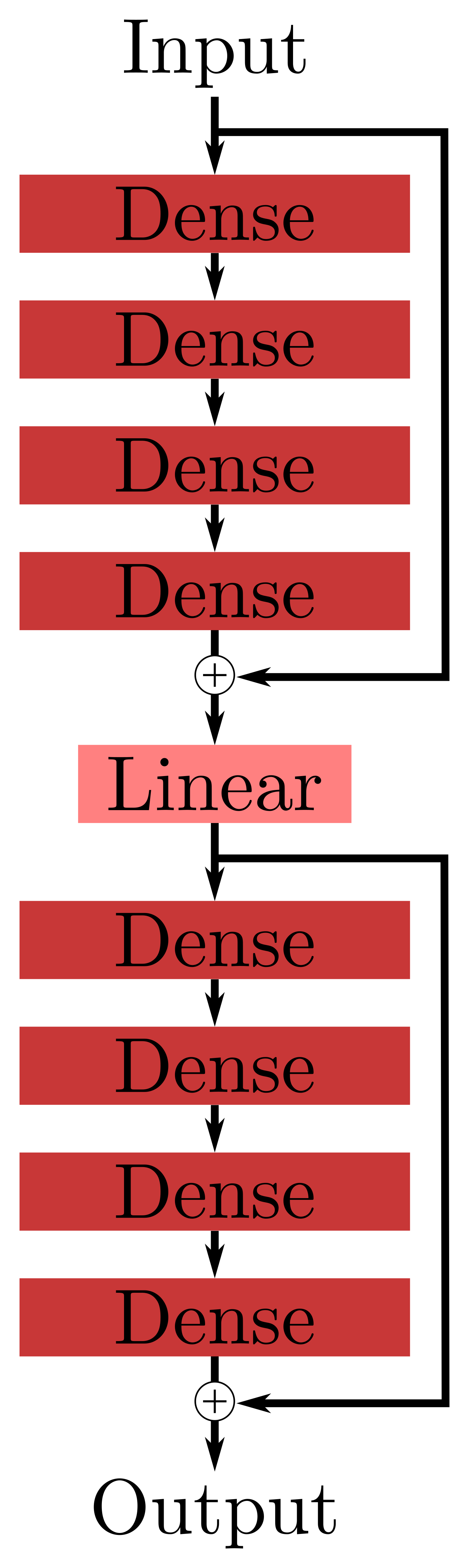}
    \vspace{-.15in}
    \caption{Layer-by-layer autoencoder architecture for 1D problems.}
    \label{fig:1d_layer_arch}
\end{figure}

\section{Neural Network Implementation Details}\label{sec:implementation}

The model training procedure is kept constant for all of the examples in this work. The networks are optimized with an Adam optimizer ($\beta_1=0.9$, $\beta_2=0.999$). Every numerical experiment starts by training a set of 20 models for a `small' number of epochs. Each of the 20 models has a randomly selected learning rate for the Adam optimizer, uniformly selected between $10^{-2}$ and $10^{-5}$. The initial training period consists of two phases: autoencoder-only (75 epochs) and full model (250 epochs). The autoencoder-only phase only enforces the autoencoder losses $\mathcal{L}_1$ and $\mathcal{L}_2$ during backpropagation. A checkpoint algorithm is used to keep track of the model with the lowest overall loss throughout the training procedure. At the end of the initial period, the best model is selected and the others are eliminated. The best model is trained for an additional 2500 epochs.

There are two network architectures in this work. 
The architectures depicted in Figures \ref{fig:2D_architecture} and \ref{fig:2D_architecture2} are applied to the two-dimensional nonlinear Poisson BVP.
The architecture depicted in Figure \ref{fig:1d_layer_arch}, is applied to one-dimensional problems. 

The two architectures have a few training variables in common. Both models use variance scaling initialization, $\ell_2$ regularization ($\lambda=10^{-6}$), and ReLu activation functions for fully connected (1D architecture) and convolutional (2D architecture) layers. Notably, the two layers immediately before and after the latent space do not have activation functions. A normalized mean squared error loss function is used for all of the loss functions, as described in Section \ref{sec:dae-bvp}. The models are trained in batches of 64 samples.

The 2D architecture utilizes convolutional layers and pooling layers, as shown in Figures \ref{fig:2D_architecture} and \ref{fig:2D_architecture2}. All convolutional layers use a kernel size of $4 \times 4$. There are differences between the convolutional layers in the encoder and the convolutional layers in the decoder. The encoder convolutional layers use a stride size of $1 \times 1$ and an increasing number of filters ($8, 16, 32, 64$). The deconvolutional layers use a stride size of $2 \times 2$ with decreasing filter size ($32, 16, 8$). Pooling layers are similar for both the encoder and decoder with a stride size of $2 \times 2$ and a pool size of $2 \times 2$.

\section{Additional Results}\label{sec:additional_results}

The repeatability of the results and models learned by the DeepGreen architecture are interesting to study from the perspective of operator convergence and latent space representations. In both cases, we aim to investigate the convergence of the model parameters to determine if the learned latent spaces and operators are unique or non-unique.

\subsection{Operator Initialization}

\begin{figure}[t]
    \centering
    \includegraphics{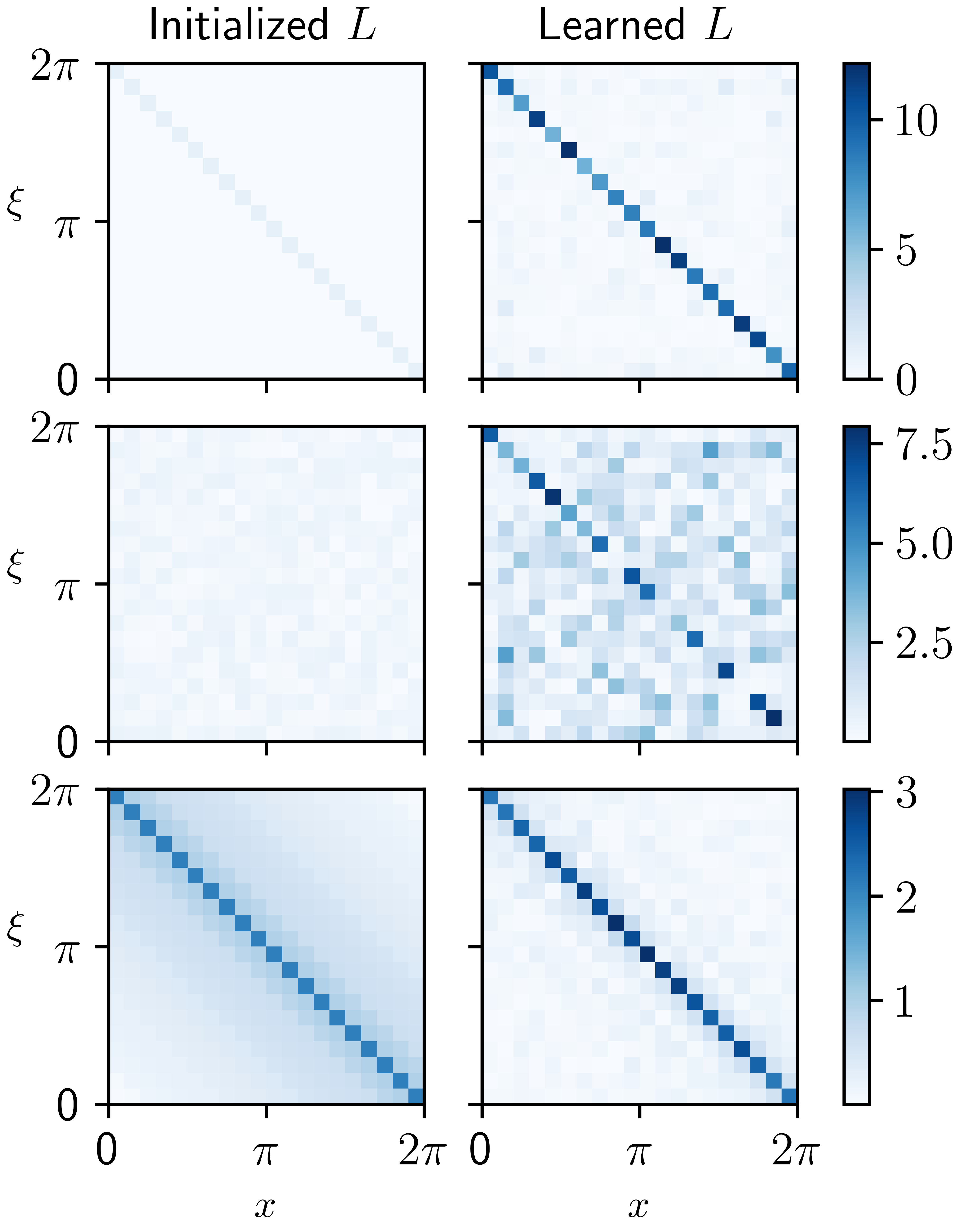}
    \vspace{-.15in}
    \caption{Initial vs learned operators for an operator matrix $L$ for different initial conditions. The top row shows identity matrix initialization, the middle row shows random initialization (He normal), and the bottom row shows a Toeplitz gradient initialization.}
    \label{fig:initial_vs_learned_ops}
\end{figure}

We repeat the training procedure for DeepGreen with three different initialization approaches for the operator $L$. Again, we train with data from the example nonlinear cubic Helmholtz model. This experiment focuses on comparing the initial values of the operator $L$ with the final values of the operator at the end of training to determine if the DeepGreen approach tends to converge to a specific operator construction. The results in Figure \ref{fig:initial_vs_learned_ops} show the initial and final operator for identity-initialized, randomly initialized, and Toeplitz-initialized operator matrices. Impressively, the result shows that the network tends to learn operators with diagonal dominance for all of the tested initialization strategies. This approach, which DeepGreen appears to prefer, draws strong parallels to the coordinate diagonalization approach commonly used in physics.

\subsection{Latent Space Analysis}

\begin{figure*}[t]
    \centering
    \includegraphics{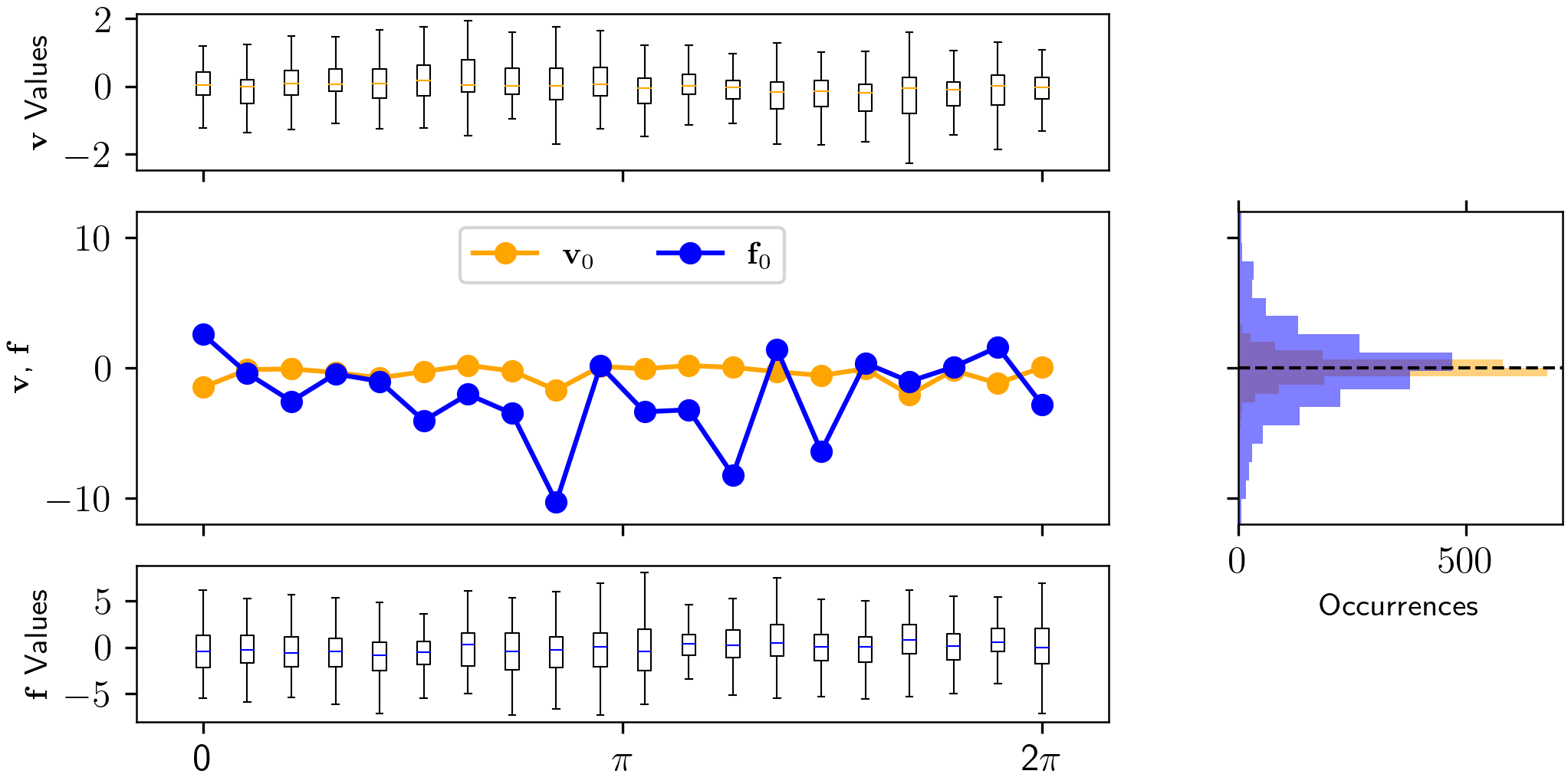}
    \caption{Statistics of latent space values of a single sample over 100 experimental runs.}
    \label{fig:latent_stats}
\end{figure*}

We repeat the training procedure for the example system, the nonlinear cubic Helmholtz model, a total of one hundred times. A single sample was selected from the training data and the latent space representation, $\mathbf{v}_i$ and $\mathbf{f}_i$, of the input vectors $\mathbf{u}_i$ and $\mathbf{F}_i$ are computed. Statistics for the latent space representations are presented in Figure \ref{fig:latent_stats}. It is evident that the latent space vectors are not identical between runs, and that the values in the vector do not follow any particular statistical distribution. This information implies that the learned weights in the model, and the learned latent space representations, vary for each training instance and do not appear to converge to a single representation.

\subsection{Residual network architecture}

All of the autoencoders used in this work use a residual network (ResNet) architecture. In order to demonstrate the advantage of the ResNet architecture, we trained six models using the DeepGreen architecture for each of the four systems. Three of the models use the ResNet skip connections, while three do not use the ResNet architecture. 

For the two simplest systems, the nonlinear cubic Helmholtz equation and the nonlinear Sturm--Liouville equation, the difference between the models with and without the ResNet skip connections was negligible. For the nonlinear cubic Helmholtz equation, the mean validation loss for the non-ResNet models was $2.7 \times 10^{-3}$ and the median validation loss was $2.4 \times 10^{-3}$. Using the ResNet architecture resulted in a mean validation loss of $3.5 \times 10^{-3}$ and a median validation loss of $8.8 \times 10^{-4}$. The ResNet architecture resulted in a lower median validation loss but a higher mean due to one of the three models performing much more poorly than the other two. The results for the nonlinear Sturm--Liouville system are analogous. With a non-ResNet architecture, the mean validation loss was $4.5 \times 10^-3$ and the median validation loss was $4.0 \times 10^{-3}$. With a ResNet architecture, the mean validation loss was $5.7 \times 10^-3$ and the median validation loss was $3.1 \times 10^{-3}$. Therefore, the use of the ResNet architecture produced similar results to a non-ResNet architecture for these two simple systems.

For the two systems that had larger losses --- the nonlinear biharmonic equation in 1D and the 2D nonlinear Poisson equation --- the ResNet architecture was clearly superior to a non-ResNet architecture. For the nonlinear biharmonic equation, the ResNet architecture yields a mean validation loss of $2.5 \times 10^{-2}$ and median validation loss of $2.8 \times 10^{-2}$ for the three models, compared with $3.8 \times 10^{-2}$ and $4.0 \times 10^{-2}$, respectively, for the non-ResNet architecture. Therefore, the ResNet architecture performed better in terms of both the mean and median loss. The ResNet architecture is absolutely vital for the nonlinear Poisson system. Without the ResNet architecture, the model essentially did not converge. Both the mean and median validation losses were $1.9 \times 10^0$. In contrast, the ResNet architecture had a mean validation loss of $1.8 \times 10^{-2}$ and a median of $1.9 \times 10^{-3}$.

\end{appendices}

\end{document}